%% file: ghyson.tex
\documentclass{article}

\usepackage[english]{babel}
\usepackage[utf8]{inputenc}
\usepackage[T1]{fontenc}
\usepackage{lmodern}

\usepackage{algorithm}
\usepackage{algpseudocode}
\usepackage{macros}

\DeclareSymbolFont{stmry}{U}{stmry}{m}{n}

\title{Computing Flowpipe of Nonlinear Hybrid Systems\\with Numerical
  Methods}

\author{Olivier Bouissou and Samuel Mimram
  \\
  {\small CEA Saclay Nano-INNOV Institut CARNOT, Gif-sur-Yvette France}
  \\[10pt]
  Alexandre Chapoutot\\
  {\small U2IS -- ENSTA ParisTech, Palaiseau France}}

\date{January 2013}

\begin{document}

\maketitle

\begin{abstract}
  \input{abstract.tex}
\end{abstract}

\section{Introduction}
\label{sec:introduction}

\input{introduction.tex}

\section{Preliminaries}
\label{sec:numerical-simulation}

\input{numerical-simulation}

\section{Guaranteed Simulation Methods}
\label{sec:guaranteed-simulation-methods}

The elaboration of our algorithm consisted essentially in adapting
simulation algorithms such as the one described in
Section~\ref{sec:numerical-simulation} in order to (i) compute with
\emph{sets of values} instead of values, and (ii) ensure that the
resulting algorithm is \emph{guaranteed} in the sense that the
set~$\hat x_k$ of values computed for $x$ at instant $t_k$ always
contains the value of the mathematical solution at instant $t_k$. This
means that we have to take in account numerical errors due to the
integration method and the use of floats (see
Section~\ref{sec:guaranteed-simulation-methods}), and design an
algorithm computing an over-approximation of jump times
(Section~\ref{sec:reachability-algorithm}). In this section, we first
briefly present our encoding of sets using affine arithmetic
(Section~\ref{sec:computing-with-sets}) and show how explicit
Runge-Kutta like numerical integration methods
(Section~\ref{sec:guaranteed-numerical-integration}) and the
polynomial interpolation
(Section~\ref{sec:guaranteed-polynomial-interpolation}) can be turned
into guaranteed algorithms.

\input{guaranteed-integration.tex}

\input{polynomial-interpolation.tex}

\input{reachability-algorithm.tex}

\section{Experimentation}
\label{sec:experimentation}

\input{experimentation}

\section{Conclusion}
\label{sec:conclusion}

\input{conclusion.tex}

\bibliographystyle{abbrv}
\bibliography{bib}

\newpage
\appendix
\section{Other examples}
Because of space constraints, we did not include the description of some
examples in the article, they can be found below.

\subsection{The Bouncing Pendulum}
This hybrid system describes pendulum attached to a rope of length $l=1.2$
falling under a gravity of $g=9.81$. The angle~$\theta$ of the pendulum (\wrt
vertical) is described by the flow equation
\[
\ddot\theta=-\frac gl\sin(\theta)
\qquad\qquad
\theta(0)=[1,1.05]
\]
The pendulum bounces on a wall when $\theta=-0.5$, in which case the reset
condition is $\dot\theta=-\dot\theta$. The guaranteed simulation of the system
produces:
\[
\includegraphics{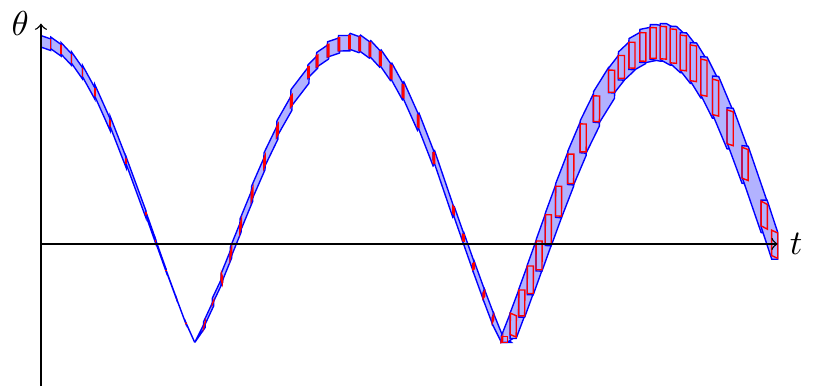}
\]
As illustration, we give here the description of the system given as input to
HySon:
{\small
\begin{verbatim}
set duration = 3.8;
set dt = 0.05;
set max_dt = 0.1;
set scope_xy = true;

init theta = [1.,1.05];
init dtheta = 0.;
init t = 0;

l = 1.2;
g = 9.81;
theta' = dtheta;
dtheta' = -g/l*sin(theta);
t' = 1;

on sin(theta) <= -0.5 do { print("Bouncing!\n"); dtheta = -dtheta };

output(t,theta);
\end{verbatim}
}

\noindent
Notice that the dynamics of the system is nonlinear (because of the presence
$\sin(\theta)$ in the flow equation) and the guard is also non linear, which
makes that it cannot be simulated with Flow$^*$.

\subsection{Wolfgram}
The simulation produced on the Wolfgram example is
\[
\includegraphics{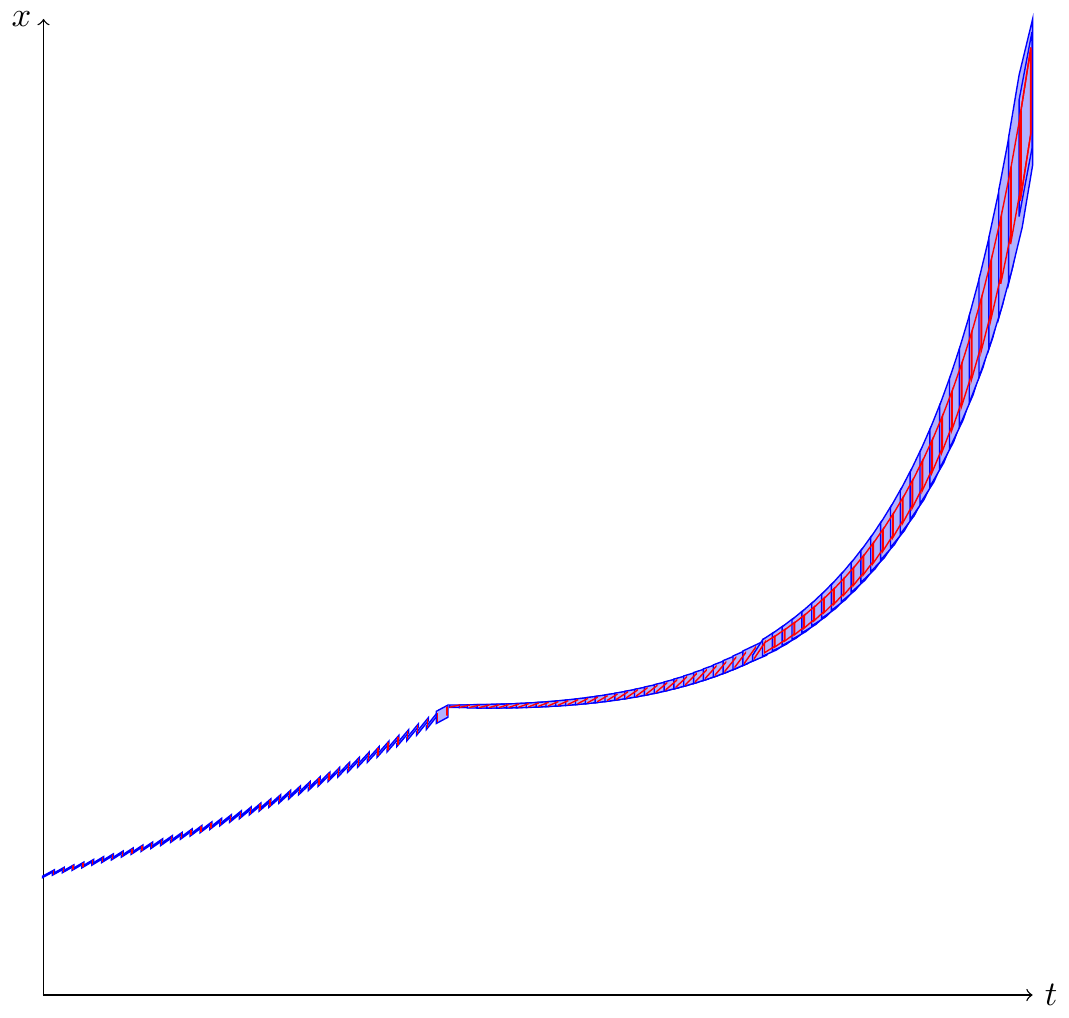}
\]

\newpage
\section{Comparison with SpaceEx}
\label{sec:comparison-spaceex}
Since the main novelty of HySon is to handle efficiently non-linear systems, we
did not detail experiments on linear ones. However, performances are comparable
with the state-of-the-art guaranteed simulators dedicated to linear systems. As
illustration, we compare here HySon with SpaceEx~\cite{spaceex} on two examples.

\subsection{Bouncing Ball}
\[
\includegraphics{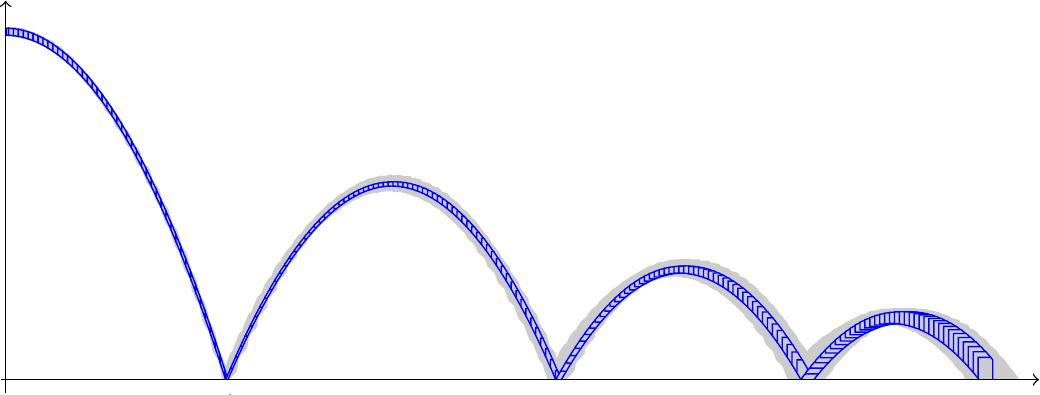}
\]
The above figure shows the flowpipe computed by SpaceEx (in gray) and by HySon
(blue polygons) for the classical bouncing-ball example, up to $t_f=20$. The
computation times were 1.031s for HySon and 1.15s for SpaceEx (we used the
support-function representation of sets using 50 directions). Notice that the
flowpipe computed by HySon is within the flowpipe of SpaceEx; we could get a
more precise results with SpaceEx by increasing the number of directions, but at
the cost of higher computation times (8.65s for 200 directions for example).

\subsection{Thermostat}
\[
\includegraphics[width=0.9\textwidth]{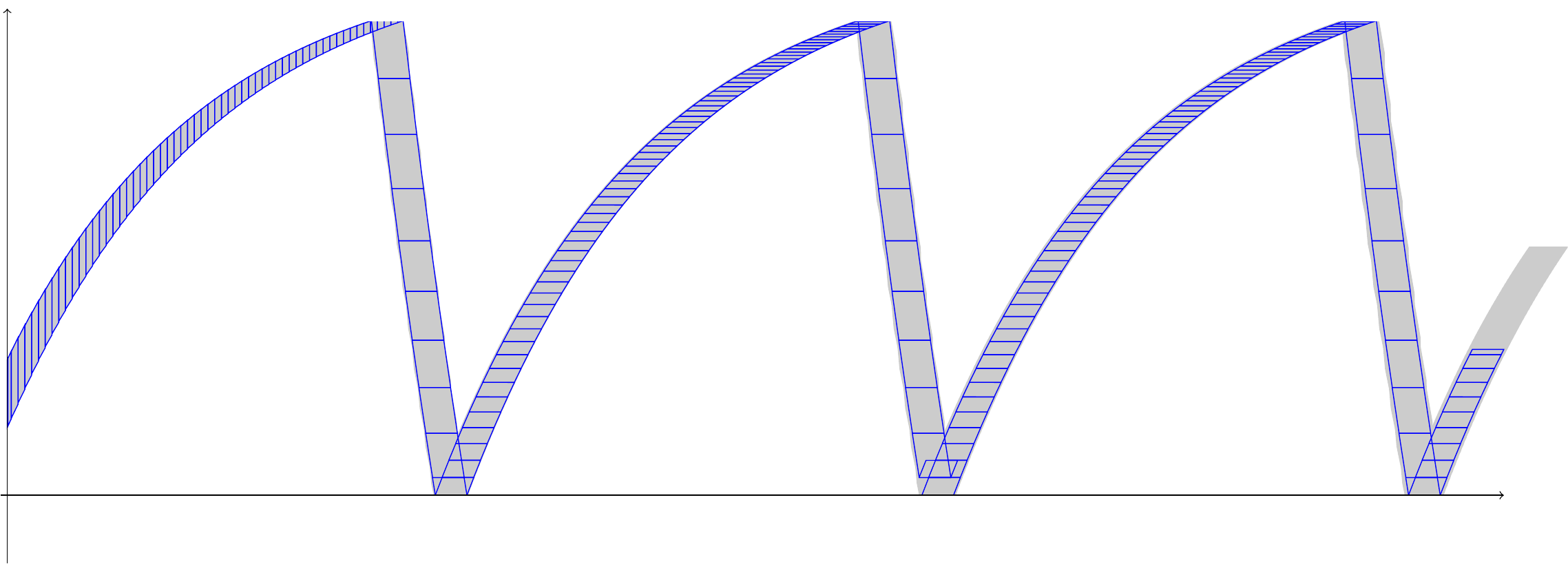}
\]
The above figure shows the flowpipe computed by SpaceEx (in gray) and by HySon
(blue sets) for the classical thermostat example, up to $t_f=15$. The
computation times were 0.89s for HySon and 0.91s for SpaceEx (we used the
support-function representation of sets using 50 directions). Notice that both
flowpipes are almost identical.
\end{document}

%% file: abstract.tex
Modern control-command systems often include controllers that perform
nonlinear computations to control a physical system, which can
typically be described by an hybrid automaton containing
high-dimensional systems of nonlinear differential equations. To prove
safety of such systems, one must compute all the reachable sets from a
given initial position, which might be uncertain (its value is not
precisely known). On linear hybrid systems, efficient and precise
techniques exist, but they fail to handle nonlinear flows or jump
conditions. In this article, we present a new tool name HySon which
computes the flowpipes of both linear and nonlinear hybrid systems
using guaranteed generalization of classical efficient numerical
simulation methods, including with variable integration step-size. In
particular, we present an algorithm for detecting discrete events
based on guaranteed interpolation polynomials that turns out to be
both precise and efficient. Illustrations of the techniques developed
in this article are given on representative examples.

%% file: introduction.tex
Modern control-command software for industrial systems are becoming
more and more complex to design. On the one side, the description of
the physical system that must be controlled, a power plant for
instance, is frequently done using partial differential equations or
nonlinear ordinary differential equations, whose number can grow very
fast when one tries to have a precise model. On the other side, the
complexity of the controller also increases when one wants it to be
precise and efficient. In particular, adaptive controllers (which
embed information on the plant dynamics) are more and more used: for
such systems, the controller may need to compute approximations of the
plant evolution using a look-up table or a simple approximation scheme
as in~\cite{bouissou_et_al_cav_09}. As an extreme example, consider a
controller for an air conditioning device in a car. In order to
correctly and pleasantly regulate the temperature in the car, the
controller takes information from the temperature of the engine but
also from the outside temperature and the sunshine on the car. Based
on these data, it acts on a cooling device, which is usually made of a
hot and a cold fluid circuit, and is thus described using usual
equations in fluid dynamics, which are given by high dimensional
nonlinear differential equations.

In an industrial context, the design of such control-command systems
is generally \emph{validated} by performing numerical simulations of a
high level description of the system using a \simulink{} like
formalism. Usually, some input scenarios are defined and numerical
simulation tools are used to observe the reaction of the system to
these inputs and check that they are in accordance with the
specifications. This methodology is widespread, because the methods
for numerical simulation used nowadays are very powerful and efficient
to approximate the behavior of complex dynamical systems, and scale
very well \wrt to both complexity and dimension. Simulation algorithms
are mainly based on two parts: algorithms to compute approximations of
the continuous evolutions of the system~\cite{SGT03}, and algorithms
to compute switching times~\cite{ZYM08}. Matlab/Simulink is the de
facto standard for the modeling and simulation of hybrid systems; we
recall its basics principles in
Section~\ref{sec:numerical-simulation-algorithms}, and refer the
reader to~\cite{BC12} for a complete formalization of its numerical
engine.

The main drawback of simulation is that it cannot give strong
guarantees on the behavior of a system, since it merely produces
approximations of it for a finite subset of the possible inputs. To
overcome this problem, verification techniques have been proposed on
slightly different models of hybrid systems. The most popular and used
technique is bounded model checking of hybrid
automata~\cite{Hen96,HR98,spaceex,Frehse05} that computes
over-approximations of the set of reachable states of a hybrid system
over a finite horizon. To apply such techniques on \simulink{}
industrial systems, one must first translate it into the hybrid
automata formalism (for example using techniques
from~\cite{AgrawalSimonKarsai04}), and then apply some simplifications
and linearizations to the model in order to obtain a linear hybrid
automaton for which the good techniques exist~\cite{spaceex}. This
process of \emph{linearization} can be performed
automatically~\cite{DT11}, but increase largely the number of discrete
states (exponentially \wrt dimension), so that we believe that it is
not applicable for large and highly nonlinear systems with stiff
dynamics.

\paragraph*{Contribution.}
In this article, we propose a new method to compute bounded horizon
over-approximations of the trajectories of hybrid systems. This method
improves our previous work~\cite{BCM12} as it modifies numerical
simulation algorithms to make them compute \emph{guaranteed} bounds of
the trajectories. Our algorithm is general enough to handle both
nonlinear continuous dynamics and nonlinear jump conditions (also
named zero-crossing events in \simulink{}). In short, our algorithm
relies on two guaranteed methods: the continuous evolution is
over-approximated using guaranteed integration of differential
equations, using a generalization of~\cite{BM06}, and the discrete
jumps are solved using a new method (presented in
Section~\ref{sec:guaranteed-polynomial-interpolation}) that can be
seen as a guaranteed version of the zero-crossing algorithm of
\simulink{}.

\paragraph*{Related work.}
We already mentioned the work on reachability analysis in hybrid
automata, either for the linear case~\cite{LGG09,spaceex}, or in the
nonlinear case where a hybridization is used to construct an
over-approximated linear automata~\cite{DT11}. Our approach is quite
different as the algorithms we propose do not suppose anything about
the differential equations and the jump conditions except their
continuity \wrt{} state space variables. Previous works also used
guaranteed numerical methods for reachability analysis of hybrid
systems~\cite{HHMW00,ERNF11}. These methods mainly use intervals as
representation of sets, such as in the library
\VNODE~\cite{Nedialkov}, to compute guaranteed bounds on the
continuous trajectories, and interval methods or a \SAT solver to
safely over-approximate the discrete jumps. Our method uses a more
expressive domain for representing sets (affine
forms~\cite{GP11,BCM12}) and polynomial interpolation for discrete
jumps, which offers an efficient bisection method.
Finally, the work closest to our is~\cite{CAS12}, in which a flowpipe
for nonlinear hybrid systems is computed using a Taylor model to
enclose the continuous behavior, and the discrete jumps are handled by
doing the intersection of elements of the Taylor model and polyhedric
guards. Compared to our approach, this work only allows for polynomial
dynamics and polyhedral guards, while we have no such restrictions (as
exemplified in Section~\ref{sec:experimentation}). Beside, as will be
clear from our benchmarks, the use of affine forms and numerical
methods is generally more efficient than Taylor models.

\paragraph*{Outline of the paper.}
The rest of this article is organized as follows. In
Section~\ref{sec:numerical-simulation}, we present our formalism for
hybrid systems and recall traditional method for their numerical
simulation. Then, in Section~\ref{sec:guaranteed-simulation-methods},
we explain how we could turn these methods into guaranteed methods
that compute enclosures rather than approximations. In
Section~\ref{sec:reachability-algorithm}, we present our main
algorithm for computing safe bounds on the trajectories of hybrid
systems, and Section~\ref{sec:experimentation} presents some
benchmarks that include both nonlinear dynamics and nonlinear jump
conditions.

%% file: numerical-simulation.tex
\subsection{Hybrid Automata}
\label{sec:state-space-formalism}

In this article, in order to facilitate the understanding of our
method, we consider hybrid systems described as hybrid
automata~(\HA). However, our tool HySon uses a slightly different
representation as in our previous work~\cite{BC12,BCM12}. This
state-space representation, comparable to the one used
in~\cite{Goebel04hybridsystems}, can encode both \HA and \simulink{}
models, as shown in~\cite{BC12}.
We denote by $\R$ the set of real numbers, and by $\B$ the set of
booleans (containing two elements, $\top$ meaning true and $\bot$
meaning false). Given a function $x:\R\to\R^n$, we denote by $x^-(t)$
its left-limit.

\begin{definition}[Hybrid automaton,~\cite{Hen96}]
  \label{def:definition-ha}
  An $n$-dimensional \emph{hybrid automaton} $\mc{H}=(L, F,E, G, R)$
  is a tuple such that $L$ is a finite set of \emph{locations}, the
  function $F:L\to(\R\times\R^n\to\R^n)$ associates a \emph{flow
    equation} to each location, $E\subseteq L\times L$ is a finite set
  of \emph{edges}, $G:E\to(\R^n\to\B)$ maps edges to \emph{guards} and
  $R:E\to(\R\times\R^n\to\R^n)$ maps edges to \emph{reset maps}.
\end{definition}

\noindent
Notice that to simplify the presentation of our approach, we consider
\HA without invariants in each location, we will discuss this point in
the conclusion. Besides, we assume that a transition $e=(l,l')$ is
taken \emph{as soon as} $G(e)$ is true.

\begin{example}
  \label{ex:windy-ball}
  We consider a modification of the classical bouncing-ball system
  that we call the \emph{windy ball}: the ball is falling but there is
  in addition an horizontal wind which varies with time. So, the
  dynamics of the horizontal position $x$ and height $y$ of the ball
  are given by
  \[
  \dot{x}(t)=10 (1 + 1.5 \sin(10t))
  \qquad\qquad
  \dot{y}(t)=v_y(t)
  \qquad\qquad
  \dot{v_y}(t)=-g
  \]
  The \HA thus has only one location $l$ such that $F(l)$ is the above
  flow. There is also one edge $e=(l,l)$ for when the ball bounces on
  the floor, with a guard $G(e)=y\leq 0$ and a reset $R(e)=
  \pa{x,y,v_y}\mapsto \pa{x,y,-0.8v_y}$.
\end{example}

The operational semantics~\cite{Hen96} of an \HA is a transition
system with two kinds of transitions for the time elapse and the
discrete jumps. From this operational semantics we can define the
trajectories of the \HA, as in~\cite{Goebel04hybridsystems}.

\begin{definition}[Trajectory of an hybrid automaton]
  \label{def:semantics-ha}
  Suppose fixed an \HA $\mc{H}=(L,F,E,G,R)$. A \emph{state} of
  $\mc{H}$ is a couple $(x,l)$ with $x\in\R^n$ and $l\in L$. A
  \emph{trajectory} of $\mc{H}$, on the time interval $[t_0,t_f]$,
  starting from an initial state $(x_0,l_0)$, is a pair of functions
  $(x,l)$ with $x:[t_0,t_f]\to\R^n$ and $l:[t_0,t_f]\to L$, such that
  there exists time instants $t_0\leq t_1\leq\ldots\leq t_n=t_f$
  satisfying, for every index~$i$,
  \begin{enumerate}
  \item $x$ is continuous and $l$ is constant on $[t_i,t_{i+1}[$,
  \item $x(0)=x_0$, $l(0)=l_0$,
  \item $\forall t\in[t_i,t_{i+1}[$, $\dot x(t)=F(l(t))(t,x(t))$,
  \item $\forall t\in]t_i,t_{i+1}[$, $\forall e=(l(t),l')\in E$,
    $G(e)(x(t))=\bot$,
  \item $G(e)(x^-(t_i))=\top$ with $e=(l^-(t_i),l(t_i))$ and
    $x(t_i)=R(e)(t_i,x^-(t_i))$.
  \end{enumerate}
\end{definition}

\noindent
In the above definition, the equations constraint the function~$x$ so
that it conforms to the flow and jump conditions
of~$\mc{H}$. Equation~2 ensures that~$x$ satisfies the initial
conditions, \eq~3 specifies that the dynamics of~$x(t)$ is the flow at
location~$l(t)$, \eq~4 and 5 ensure that the $t_i$ are the instants
where jumping conditions occur and that $x$ evolves as described by
reset maps when the corresponding guard is satisfied. Notice that we
do not consider Zeno phenomena here as we assume that there are
finitely many jumps between $t_0$ and $t_f$. Also, we do not discuss
conditions ensuring existence and unicity of trajectories as this is
beyond the scope of this paper~\cite{HNW93}, but implicitly suppose
that these are granted. We suppose fixed initial and terminal
simulation times $t_0$ and $t_f$. Given an \HA $\mc{H}$ and an initial
state $(x_0,l_0)$, we denote by $Reach_{\mc{H}}(x_0,l_0)$ the
continuous trajectory on $[t_0,t_f]$ as defined above, and given
$X_0\subseteq\R^n$ and $L_0\subseteq L$, we define
$Reach_{\mc{H}}(S_0,L_0)=\bigcup_{x_0\in X_0,l_0\in
  L_0}Reach_{\mc{H}}(x_0,l_0)$.

Computing the set $Reach_{\mc{H}}(X_0,L_0)$ for an \HA $\mc{H}$ is
sufficient in order to decide the reachability of some region in the
state space, and thus often to prove its safety (for bounded time). As
trajectories are in general not computable, over-approximations must
be performed: this is the goal of our algorithm presented in
Sections~\ref{sec:guaranteed-simulation-methods}
and~\ref{sec:reachability-algorithm}. In
Section~\ref{sec:numerical-simulation-algorithms}, we present
numerical algorithms, used for example by \simulink{}, that allow to
compute \emph{approximations} of the set $Reach_{\mc{H}}(x_0,l_0)$ for
some initial state $(x_0,l_0)\in\R^n\times L$. In
Section~\ref{sec:guaranteed-simulation-methods} we present how we can
adapt these methods in order to be safe \wrt the exact trajectories
of~$\mc{H}$.

\subsection{Numerical Simulation}
\label{sec:numerical-simulation-algorithms}

Numerical simulation aims at producing discrete approximations of the
trajectories of an hybrid system $\mc{H}$ on the time interval
$[t_0,t_f]$. We described in details in~\cite{BC12} how the simulation
engine of \simulink{} operates, and briefly adapt here this simulation
engine to \HA.

Suppose that $\mc{H}$ is an \HA, $(x_0,l_0)$ an initial state, and
$(x(t),l(t))$ a trajectory of $\mc{H}$ starting from $(x_0,l_0)$. A
numerical simulation algorithm computes a sequence
$(t_k,x_k,l_k)_{k\in[0,N]}$ of time instants, variables values and
locations such that $\forall k\in[0,N],\ x_k\approx x(t_k)$. Most of
the difficulty lies in approximating the discrete jumps (instants
where a guard becomes true), which are called \emph{zero-crossings} in
the numerical simulation community. In order to compute
$(t_k,x_k,l_k)$, the following simulation loop is used, where $h_k$ is
the simulation step-size (that can be modified to a smaller value in
order to maintain a good precision): \vbox{
\begin{algorithmic}[1]
  \Repeat
  \Let{$x_{k+1}$}{$\text{SolveODE}(F(l(t_k),x_k,h_k)$}\Comment{Solver
    step 1}
  \Let{$(x_{k+1},l_{k+1})$}{$\text{SolveZC}(x_k,x_{k+1})$}\Comment{Solver
    step 2} \State{\text{compute} $h_{k+1}$} \Let{$k$}{$k+1$}
  \Until{$t_k\ge t_f$}
\end{algorithmic}
}
\noindent
In this simulation loop, the solver first makes a continuous
transition between instants $t_k$ and $t_k+h_k$ under the assumption
that no jump occurs (solver step 1), and then it verifies this
assumption (solver step 2). If it turns out that there was a jump
between $t_k$ and $t_k+h_k$, the solver approximates as precisely as
possible the time $t\in[t_k,t_k+h_k]$ at which this jump occurred. We
briefly detail both steps in the rest of this section.

\paragraph*{Solver step 1.}
The continuous evolution of $x$ between $t_k$ and $t_k+h_k$ is
described by $\dot{x}(t)=F(l(t_k))\bigl(t,x(t)\bigr)$ and
$x(t_k)=x_k$. So, we want to compute an approximation of the solution
at $t_k+h_k$ of the initial value problem (\IVP), with $f=F(l(t_k))$:
\begin{equation}
  \label{eq:ivp}
  \dot{x}(t)=f(t,x(t))
  \qquad\qquad\qquad
  x(t_k)=x_k
\end{equation}
(we assume classical hypotheses on $f$ ensuring existence and
uniqueness of a solution of \IVP). Usually, precise simulation
algorithms often rely on a variable step solver, for which $(h_k)$ is
not constant. The simplest is probably the Bogacki-Shampine
method~\cite{SGT03}, also named \ODE[$23$]. It computes $x_{k+1}$ by
\begin{subequations}
  \label{eq:bogacki-shampine-ode23}
  \small
  \begin{align}
    &
    k_1=f(t_k, x_k)
    \qquad
    k_2=f(t_k + \frac{h_k}{2}, x_k + \frac{h_k}{2}k_1)
    \qquad
    k_3=f(t_k + \frac{3h_k}{4}, x_k + \frac{3h_k}{4}k_2)
    \label{eq:bogacki-shampine-ode3-start}
    \\
    &
    x_{k+1}=x_k + \frac{h_k}{9}\pa{2k_1 + 3k_2 +4k_3}
    \label{eq:bogacki-shampine-ode3-end}
    \\
    &
    k_4=f(t_k + h_k, x_{k+1})
    \qquad\qquad\qquad
    z_{k+1}=x_k+\frac{h_k}{24}\pa{7k_1+6k_2+8k_3+3k_4}
    \label{eq:bogacki-shampine-ode3-err}
  \end{align}
\end{subequations}
The value $z_{k+1}$ defined in~\eqref{eq:bogacki-shampine-ode3-err} is
a third order approximation of $x(t_k+h_k)$, whereas $x_{k+1}$ is a
second order approximation of this value, and is used to estimate the
error $\mathrm{err}=|x_{k+1}-z_{k+1}|$. This error is compared to a
given \emph{tolerance} $\mathrm{tol}$ and the step-size is changed
accordingly: if the error is smaller then the step is \emph{validated}
and the step-size increased in order to speed up computations (in
\ODE[$23$], next step-size is computed with $h_{k+1}=h_k
\sqrt[3]{\mathrm{tol}/\mathrm{err}}$), if the error is greater then
the step is \emph{rejected} and the computation is tried again with
the smaller step-size $h_k/2$. We refer to \cite[p.~167]{HNW93} for a
complete description on such numerical methods.

\paragraph*{Solver step 2.}
Once $x_k$ and $x_{k+1}$ computed, the solver checks if there were a
jump in the time interval $[t_k,t_{k+1}]$. In order to do so, it tests
for each edge $e$ starting from~$l_k$ whether $G(e)(x_k)$ is false and
$G(e)(x_{k+1})$ is true. If there is no such edge, then it is
considered that no jump occurred, we set $l_{k+1}=l_k$ and continue
the simulation. Notice this technique does not guarantee the detection
of all events occurring between $[t_k,t_{k+1}]$ as explained
in~\cite{ZYM08} or~\cite{EKP01}.

If the solver finds such an edge, this means that there was a jump on
$[t_k,t_{k+1}]$ and we must approximate the first time instant
$\xi\in[t_k,t_k+h_k]$ such that $G(e)(x(\xi))$ is true. To do so, the
solver encloses $\xi$ in an interval $[t_l,t_r]$ starting with
$t_l=t_k$ and $t_r=t_k+h_k$, and reduces this interval until the time
precision $|t_l-t_r|$ is smaller than some parameter. To reduce the
width of the interval, the solver first makes a guess for $\xi$ using
a linear extrapolation and then computes an approximation of $x(\xi)$
using a polynomial interpolation of $x$ on $[t_k,t_k+h_k]$. Depending
on $G(e)(x(\xi))$, it then sets $t_l=\xi$ or $t_r=\xi$ and starts
over. In the case of Hermite interpolation (which is the method used
together with the \ODE[$23$] solver), the polynomial interpolation is
given, for $t\in[t_k,t_k+h_k]$, by
\begin{equation}
  \label{eq:hermite}
  x(t)\approx (2\tau^3-3\tau^2+1)x_k+(\tau^3-2\tau^2+\tau)h_k\dot{x}_k+(-2\tau^3+3\tau^2)x_{k+1}+(\tau^3-\tau^2)h_k\dot{x}_{k+1}
\end{equation}
where $\tau=(t-t_k)/h_k$, and $\dot{x}_k$, $\dot{x}_{k+1}$ are
approximations of the derivative of $x$ at $t_k$, $t_{k+1}$. For more
details on zero-crossing algorithms, we refer to~\cite{BC12,ZYM08}.

\begin{example}
  Consider the windy ball again (Example~\ref{ex:windy-ball}). The red
  curve below is the result of the simulation for $t\in[0,13]$ using
  \simulink. In blue is the flowpipe computed by \hyson whose
  computation is going to be described in next sections.
  \[
  \includegraphics{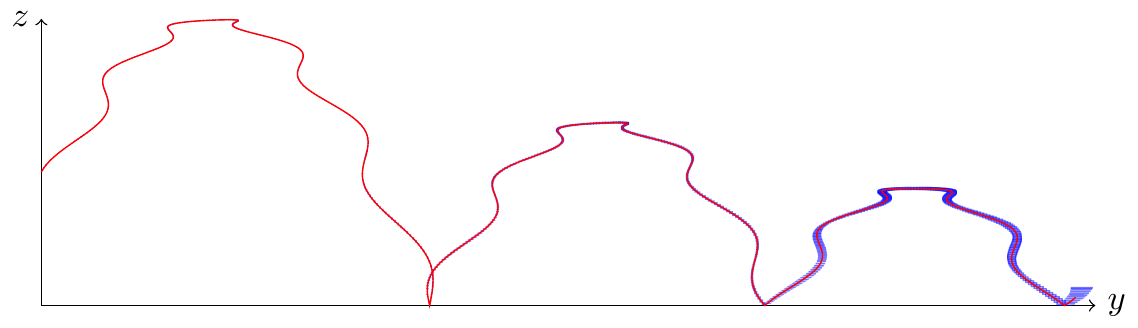}
  \]
\end{example}

%% file: guaranteed-integration.tex
\subsection{Computing with Sets}
\label{sec:computing-with-sets}
The simplest and most common way to represent and manipulate sets of
values is \emph{interval arithmetic}~\cite{Moore66}. Nevertheless,
this representation usually produces too much over-approximated
results, because it cannot take dependencies between variables in
account: for instance, if $x=[0,1]$, then $x-x=[-1,1]\neq0$. More
generally, it can be shown for most integration schemes that the width
of the result can only grow if we interpret sets of values as
intervals.

To avoid this problem we use an improvement over interval arithmetic
named \textit{affine arithmetic}~\cite{AffineA97} which can track
linear correlations between program variables. A set of values in this
domain is represented by an \textit{affine form}~$\hat x$ (also called
a \textit{zonotope}), which is a formal expression of the form $\hat
x=\alpha_0+\sum_{i=1}^n \alpha_i\varepsilon_i$ where the coefficients
$\alpha_i$ are real numbers, $\alpha_0$ being called the \emph{center}
of the affine form, and the $\varepsilon_i$ are formal variables
ranging over the interval $[-1,1]$.  Obviously, an
interval~$a=[a_1,a_2]$ can be seen as the affine form $\hat
x=\alpha_0+\alpha_1\varepsilon$ with $\alpha_0=(a_1+a_2)/2$ and
$\alpha_1=(a_2-a_1)/2$. Moreover, affine forms encode linear
dependencies between variables: if $x\in[a_1,a_2]$ and $y$ is such
that $y=2x$, then $x$ will be represented by the affine form~$\hat x$
above and $y$ will be represented as $\hat
y=2\alpha_0+2\alpha_1\varepsilon$.

Usual operations on real numbers extend to affine arithmetic in the
expected way. For instance, if
\begin{math}
  \hat x=\alpha_0+\sum_{i=1}^n\alpha_i\varepsilon_i
\end{math}
and
\begin{math}
  \hat y=\beta_0+\sum_{i=1}^n\beta_i\varepsilon_i
\end{math}, then with $a,b,c\in\R$ we have
\begin{math}
  a \hat x+ b \hat y + c =
  (a\alpha_0+b\beta_0+c)+\sum_{i=1}^n(a\alpha_i+b\beta_i)\varepsilon_i
\end{math}.  However, unlike the addition, most operations create new
noise symbols. Multiplication for example is defined by
\begin{math}
  \hat x\times\hat y = \alpha_0\alpha_1+\sum_{i=1}^n (\alpha_i\beta_0
  + \alpha_0 \beta_i)\varepsilon_i + \nu \varepsilon_{n+1}
\end{math}, where $\nu=\left(\sum_{i=1}^n |\alpha_i|\right)\times
\left(\sum_{i=1}^n|\beta_i|\right)$ over-approximates the error
between the linear approximation of multiplication and multiplication
itself. Other operations, like $\sin$, $\exp$, are evaluated using
their Taylor expansions. The set-based evaluation of an expression
only consists in interpreting all the mathematical operators (such as
$+$ or $\sin$) by their counterpart in affine arithmetic. We will
denote by $\Aff(e)$ the evaluation of the expression $e$ using affine
arithmetic, see \cite{BCM12} for practical implementation details.

\subsection{Guaranteed Numerical Integration}
\label{sec:guaranteed-numerical-integration}

Recall from Section~\ref{sec:numerical-simulation} that a numerical
integration method computes a sequence of approximations $(t_n,x_n)$
of the solution $x(t;x_0)$ of the \IVP defined in~\eqref{eq:ivp} such
that $x_n \approx x(t_n;x_0)$. Every numerical method member of the
Runge-Kutta family follows the \textit{condition
  order}~\cite[Chap.~II.2, Thm.~2.13]{HNW93}. This condition states
that a method is of order $p$ if and only if the $p+1$ first
coefficients of the Taylor expansion of the true solution and the
Taylor expansion of the numerical method are equal. The
\emph{truncation error} measures the distance between the true
solution and the numerical solution and it is defined by $x(t_n;x_0) -
x_n$. Using the condition order, it can be shown that this truncation
error is proportional to the Lagrange remainders. We now briefly
recall our approach to make any explicit Runge-Kutta method
guaranteed, which is based on this observation, see~\cite{BCD13} for a
detailed presentation.

The general form of an explicit \emph{$s$-stage Runge-Kutta formula}
(using $s$ evaluations of $f$) is
\[
x_{n+1} = x_n + h \sum_{i=1}^s b_i k_i \qquad \text{with} \qquad k_i =
f\Bigl(t_n + c_i h,\, x_n + h \sum_{j=1}^{i-1} a_{ij}k_j\Bigr)
\]
for $1\leq i\leq s$. The coefficients $c_i$, $a_{ij}$ and $b_i$ are
usually summarized in a Butcher table (see~\cite{HNW93}) which fully
characterizes a Runge-Kutta method. We denote by
$\phi(t)=x_n+h_t\sum_{i=1}^sb_ik_i(t)$, where $k_i(t)$ is defined as
previously with $h$ replaced by $h_t=t-t_n$. Hence the truncation
error is defined by
\begin{equation}
  \label{eq:truncation-error}
  x(t_n;x_0) - x_n\qeq\frac{h_n^{p+1}}{(p+1)!} \left( f^{(p)}\pa{\xi,x(\xi)\right) - \frac{\diff^{p+1}\phi}{\diff t^{p+1}}(\eta)}
\end{equation}
for some $\xi\in]t_k, t_{k+1}[$ and $\eta\in]t_n,
t_{n+1}[$. In~\eqref{eq:truncation-error}, $f^{(p)}$ stands for the
$p$-th derivative of function $f$ \wrt time~$t$, and $h_n=t_{n+1}-t_n$
is the step size. In~\eqref{eq:truncation-error}, the Lagrange
remainder of the exact solution is $f^{(p)}\left(\xi,
  x(\xi;x_0)\right)$ and the Lagrange remainder of the numerical
solution is $\frac{d^{p+1}\phi}{dt^{p+1}}(\eta)$.

The challenge to make Runge-Kutta integration schemes safe \wrt the
exact solution of \IVP amounts to bounding the result
of~\eqref{eq:truncation-error}. The remainder
$\frac{d^{p+1}\phi}{dt^{p+1}}(\eta)$ is straightforward to bound
because the function $\phi$ only depends on the value of the step size
$h$, and so does its $(p+1)$-th derivative:
\begin{equation}
  \label{eq:bound-method}
  \frac{\diff^{p+1}\phi}{\diff t^{p+1}}(\eta)
  \in
  \Aff\pa{\frac{\diff^{p+1}\phi}{\diff t^{p+1}}([t_n, t_{n+1}])}
\end{equation}

However, the expression $f^{(p)}\pa{\xi, x(\xi;x_0)}$ is not so easy
to bound as it requires to evaluate $f$ for a particular value of the
\IVP solution $x(\xi;x_0)$ at a unknown time $\xi\in
]t_n,t_{n+1}[$. The solution we used is similar to the one found
in~\cite{Nedialkov,BM06}: we first compute an a priori enclosure of
the \IVP on the interval $[t_n,t_{n+1}]$. To do so, we use the Banach
fixed-point theorem on the Picard-Lindelöf operator~$P$, defined by
$P(x,t_n,x_n)=t\mapsto x_n+\int_{t_n}^tf(s,x(s))\diff s$. Notice that
this operator is the integral form of~\eqref{eq:ivp}, so a fixpoint of
this operator is also a solution of~\eqref{eq:ivp}.

Now, to get an a priori enclosure of the solution over
$[t_n,t_{n+1}]$, we prove that the operator $P$ (which is an operator
on functions) is contracting and use Banach theorem to deduce that it
has a fixpoint. To find the enclosure $\hat z$ on the solution, we
thus iteratively solve using affine arithmetic the equation $P(\hat
z,t_n, x_n)([t_n, t_{n+1}]) \subseteq \hat z$. Then, we know that the
set of functions $[t_n,t_{n+1}]\to \hat z$ contains the solution of
the \IVP, so $\hat z$ can be used an enclosure of the solution of \IVP
over the time interval $[t_n, t_{n+1}]$. We can hence bound the
Lagrange remainder of the true solution with $\hat z$ such that
\begin{equation}
  \label{eq:bound-true-solution}
  f^{(p)}\left(\xi, x(\xi;x_0)\right)
  \in
  \Aff\pa{f^{(p)}\left([t_n, t_{n+1}], \hat z\right)}
\end{equation}

Finally, using~\eqref{eq:bound-method}
and~\eqref{eq:bound-true-solution} we can prove
Theorem~\ref{thm:guaranteed-integration} and thus bound the distance
between the approximations point of any explicit Runge-Kutta method
and any solution of the
\IVP.
\begin{theorem}
  \label{thm:guaranteed-integration}
  Suppose that $\Phi$ is a numerical integration scheme and
  $\Phi_\Aff$ is the evaluation of $\Phi$ using affine
  arithmetic. Given a set $S_0\subseteq\R^n$ of initial states, and an
  affine form $\hat x_0$ such that $S_0\subseteq\hat x_0$, let
  $(t_n,\hat x_n)$ be a sequence of time instants and affine forms
  defined by $\hat x_{n+1}=\hat x'_{n+1}+\hat e_{n+1}$ where
  $(t_{n+1},\hat x'_{n+1})=\Phi_{\Aff}(t_n,\hat x_n)$ and~$\hat
  e_{n+1}$ is the truncation error as defined
  by~\eqref{eq:truncation-error} and is evaluated
  using~\eqref{eq:bound-method}
  and~\eqref{eq:bound-true-solution}. Then, for any $x\in S_0$ and
  $n\in\N$ we have $x(t_n;x)\in\hat x_n$.
\end{theorem}

%% file: polynomial-interpolation.tex
\subsection{Guaranteed Polynomial Interpolation}
\label{sec:guaranteed-polynomial-interpolation}
From two (guaranteed) solutions $x_n,x_{n+1}$ at times $t_n,t_{n+1}$
of an \IVP, one would like to deduce by interpolation all the
solutions $x(t)$ with $t\in[t_n,t_{n+1}]$. This question has motivated
a series of work on polynomial approximations of solutions,
a.k.a. \textit{continuous extension}, see \cite[Chap.~6]{HNW93}. We
briefly recall the polynomial interpolation method based on
Hermite-Birkhoff which is the main method used for continuous
extension. Furthermore, we present a new extension of this method
allowing us to compute a guaranteed polynomial interpolation using the
result of the Picard-Lindelöf operator.

Suppose given a sequence $(t_i, x^{(k)}_i)$ of $n+1$ computed values
of the solution of an \IVP and its derivative at instants~$t_i$, with
$0\leq i\leq n$ and $k={0,1}$. Remark that these values are those
produced by numerical integration methods.  The goal of
Hermite-Birkhoff polynomial interpolation is to build a polynomial
function
\begin{math}
  p(t) = \sum_{i=0}^n \pa{x_i A_i(t) + x^{(1)}_iB_i(t)}
\end{math}
of degree $N=2n+1$ from these values such that
\begin{math}
  A_i(t)=\pa{1 - 2(t - t_i)\ell_i'(t_i)}\ell_i^2(t)
\end{math},
\begin{math}
  B_i(t)=(t - t_i)\ell_i^2(t)
\end{math},
\begin{math}
  \ell_i(t)=\prod_{j=0,j\neq i}^n\frac{t - t_j}{t_i - t_j}
\end{math}, and
\begin{math}
  \ell'_i(t_i)= \sum_{k=0,k\neq i}^n\frac{1}{t_i - t_k}
\end{math}: the functions $\ell_i(t)$ are the Lagrange polynomials and
this interpolation generalizes the Lagrange interpolation. Under the
assumption that all the $t_i$ are distinct, we know that the
polynomial interpolation is unique. For instance,
\eq~\eqref{eq:hermite} is associated to the Hermite-Birkhoff
polynomial with $n=1$. We know that interpolation error $x(t;x_0) -
p(t)$ is defined by
\begin{math}
  \frac{x^{(N+1)}(\xi)}{(N+1)!}\prod_{i=0}^n (t-t_i)^2
\end{math} with $\xi\in [t_0, t_n]$, which can be reformulated as
\begin{equation*}
  x(t;x_0) - p(t) \qeq \frac{f^{(N)}(\xi,x(\xi))}{(N+1)!}\prod_{i=0}^n (t-t_i)^2
  \qquad\text{with $\xi\in [t_k, t_{k+1}]$}
\end{equation*}
In consequence, to guarantee the polynomial interpolation, it is
enough to know an enclosure of the solution $x(t)$ of \IVP on the
interval $[t_k, t_{k+1}]$. And fortunately, we can reuse the result of
the Picard-Lindelöf operator in that context. In next section, this
guaranteed polynomial interpolation will be used to approximate the
solution of an \IVP in order to compute jump times.

\begin{theorem}
  Let $p_{\Aff}(t)$ be the interpolation polynomial based on $n+1$
  guaranteed solutions $\hat x_i$ of an~\IVP~\eqref{eq:ivp} and $n+1$
  evaluations $\hat x_i^{(1)}$ of $f$ with affine arithmetic, and
  let~$\hat z$ be the result of the Picard-Lindelöf operator. We have,
  \[
  \forall t \in [t_k, t_{k+1}],\quad x(t)
  \in \Aff\pa{p_{\Aff}(t) +\frac{f^{(N)}([t_k, t_{k+1}],\hat z)}{(N+1)!} \prod_{i=0}^n (t-t_i)^2}
  \]
\end{theorem}

%% file: reachability-algorithm.tex
\section{Reachability Algorithm}
\label{sec:reachability-algorithm}

We present in this section our main algorithm to compute an
over-approximation of the set of reachable states of linear or
nonlinear hybrid systems
(Algorithm~\ref{alg:guaranteed-simulation-loop}), which is based on
the guaranteed numerical methods presented in
Section~\ref{sec:guaranteed-simulation-methods}. In a nutshell, it
works as follows. It produces a sequence of values $(\hat t_n,\hat
x_n^h,\hat x_n,l_n)$ such that $l_n$ is the current location, $\hat
t_n$ is a time interval, $\hat x_n$ is an over-approximation of $x(t)$
for every $t\in\hat t_n$, and $\hat x_n^h$ is an over-approximation of
$x(t)$ for every $t\in[\hat t_n,\hat t_{n+1}]$, \ie an
over-approximation of the trajectory between two discrete instants
(here $[\hat t_n,\hat t_{n+1}]$ designates the convex hull of the
union of the two affine forms $\hat t_n$ and $\hat t_{n+1}$). Our
method uses the guaranteed \ODE solver described in
Section~\ref{sec:guaranteed-numerical-integration} to compute $\hat
x_{n+1}$ and $\hat x_n^h$, and the guaranteed polynomial interpolation
of Section~\ref{sec:guaranteed-polynomial-interpolation} to precisely
and safely enclose the potential jumping times between $t_n$ and
$t_{n+1}$, and thus refine $t_{n+1}$ and $\hat x_{n+1}$.

\bigskip\noindent\textbf{Trivalent Logic.}
First, notice that since we are working with sets of values, the
evaluation of a boolean condition, such as $x\geq 0$, is not
necessarily false or true, but can also be false for some elements and
true for some other elements in the set~$\hat x$ (for instance when
$\hat x=[-1,1]$ in the preceding example). In order to take this in
account, boolean conditions are evaluated in the domain of
\emph{trivalent logic} instead of usual booleans~$\B$. This logic is
the natural extension of boolean algebra to the three following
values: $\bot$ (false), $\top$ (true) and $\botop$ (unknown). We
denote this set by $\B^*$. Notice that a function $g:\R^n\to\B$
naturally extends to a function $\Aff(g):\mc{P}(\R^n)\to\B^*$ using
affine arithmetic and trivalent logic. In particular, the guards of
the discrete jumps will be evaluated in $\B^*$, which brings
subtleties in the zero-crossing detection algorithm (when such a guard
evaluates to $\botop$), as we will see in next section. In the
following, we shall write $g$ for $\Aff(g)$ when it is clear from the
context.

\bigskip\noindent\textbf{Main Algorithm.}
Let $\mc{H}$ be an \HA as defined in
Definition~\ref{def:definition-ha}. Our method computes a sequence of
values $(\hat t_n,\hat x_n,\hat x_n^h,l_n)$ such that~$l_n$ is the
current mode of the \HA, $t_n$ is a time interval and $\hat x_n$ and
$\hat x_n^h$ are affine forms such that we have
\[
\forall t\in\hat t_n\ x(t)\in\hat x_n
\quad
\forall t\in[\hat t_n,\hat t_{n+1}],
\qquad
x(t)\in\hat x_n^h
\]
for all trajectories of $\mc{H}$. To compute this sequence, we start
from $\hat t_0=0$ and iterate until the lower bound of $\hat t_n$
(denoted $\inf(\hat t_n)$) is lower than $t_f$. The guaranteed
simulation loop is given in
Algorithm~\ref{alg:guaranteed-simulation-loop}, where $\gsolveODE{}$
is the guaranteed solver of \ODE presented in
Section~\ref{sec:guaranteed-numerical-integration} and $\gsolveZC{}$
is the procedure described below. Notice that the function
$\gsolveODE{}$ outputs both $\hat x_{n+1}$, the tight
over-approximation of $x$ at $\hat t_n+h_n$, and $\hat x_n^h$, the
result of Picard iteration (see
Section~\ref{sec:guaranteed-numerical-integration}) since we reuse it
in $\gsolveZC{}$.

\begin{algorithm}[tbp]
  \caption{Guaranteed simulation algorithm}
    \label{alg:guaranteed-simulation-loop}
  \begin{algorithmic}[1]
    \Require{$\mc{H}=(L,F,E,G,R)$, a, hybrid automaton}
    \Require{$\hat x_0$, $l_0$, $h_0$, $t_f$}\Comment{Initial state, step-size and final time}
    \Let{$n$}{$0$}
    \Let{$\hat t_n$}{$0$}
    \While{$\inf(\hat t_n)\leq t_f$}
      \Let{$(\hat x_{n+1},\hat x_n^h)$}{$\gsolveODE{}(F(l_n),\hat
        x_n,h_n)$}
      \Let{$(\hat x_{n+1},\hat t_{n+1},l_{n+1})$}{$\gsolveZC{}(l_n,\hat
        x_n,\hat x_{n+1},\hat x_n^h,\hat t_n,h_n)$}
      \Let{$n$}{$n+1$}
      \EndWhile
  \end{algorithmic}
\end{algorithm}

\bigskip\noindent\textbf{Detecting Jumps.}
We now present our algorithm ($\gsolveZC{}$) for detecting and
handling discrete jumps. Let $\mc{H}=(L,F,E,G,R)$ be an \HA, and let
$l_n$, $\hat x_n$, $\hat x_{n+1}$ and $\hat x_n^h$ be the states
computed with $\gsolveODE{}$.  Let us denote $l_n^\bullet$ the set of
all transitions originating from $l_n$, \ie $l_n^\bullet=\{e\in
E\,|\,\exists l\in L,\ e=(l_n,l)\}$.  A transition $e\in l_n^\bullet$
was \emph{surely} activated between $t_n$ and $t_n+h_n$ if $G(e)(\hat
x_n)=\bot$ and $G(e)(\hat x_{n+1})=\top$. The transition $e$ was
\emph{maybe} activated if if $G(e)(\hat x_n)=\bot$ and $G(e)(\hat
x_{n+1})=\botop$. Note that in both cases we have $G(e)(\hat
x_n^h)=\botop$. In this section, we present our algorithm in the
simple (but most common) case where we have only one transition
activated at a given time, and where we are not in the situation of
$G(e)(\hat x_n)=G(e)(\hat x_{n+1})=\bot$ with $G(e)(\hat
x_n^h)=\botop$; we discuss these two cases later.

\begin{figure}[tbp]
  \centering
  \[
  \begin{tikzpicture}[yscale=0.8]
    \draw[color=lightgray,fill=lightgray] (.5,1.4) .. controls (1.5,.1)
    .. (2,1.1) -- (2,.4) .. controls (1.5,-.6) .. (.5,.7);
    \draw[color=darkgray,fill=darkgray] (.5,1.3) .. controls (1.5,0)
    .. (2,1) -- (2,.5) .. controls (1.5,-.5) .. (.5,1);
    \draw[->] (-.5,0) -- (2.5,0);
    \draw (1,-1) node {Case a};
    \draw[->] (0,-.5) -- (0,1.5);
    \draw (.5,-.1) -- (.5,.1);
    \draw (0,1.5) node[left]{$x$};
    \draw (2.5,0) node[right]{$t$};
    \draw (.5,0) node[below]{$t_n$};
    \draw (2,-.1) -- (2,.1);
    \draw (2,0) node[below]{$t_{n+1}$};
  \end{tikzpicture}
  \qquad
  \begin{tikzpicture}[yscale=0.8]
    \draw[color=lightgray,fill=lightgray] (.5,1.4) .. controls (1.5,.1)
    .. (2,1.1) -- (2,.4) .. controls (1.5,-.6) .. (.5,.7);
    \draw[color=darkgray,fill=darkgray] (.5,1.3) .. controls (1.5,0)
    .. (2,1) -- (2,.5) .. controls (1.5,-.1) .. (.5,1);
    \draw[->] (-.5,0) -- (2.5,0);
    \draw (1,-1) node {Case b};
    \draw[->] (0,-.5) -- (0,1.5);
    \draw (.5,-.1) --(.5,.1);
    \draw (0,1.5) node[left]{$x$}; \draw (2.5,0)
    node[right]{$t$}; \draw (.5,0) node[below]{$t_n$}; \draw (2,-.1) --
    (2,.1); \draw (2,0) node[below]{$t_{n+1}$};
  \end{tikzpicture}
  \qquad
  \begin{tikzpicture}[yscale=0.8]
    \draw[color=lightgray,fill=lightgray] (.5,.8) .. controls (1.5,.3) .. (2,-.2) -- (2,-.9) .. controls (1.5,-.1) .. (.5,.2) -- cycle;
    \draw[color=darkgray,fill=darkgray] (.5,.6) .. controls (1.5,.2) .. (2,-.4) -- (2,-.6) .. controls (1.5,-.1) .. (.5,.4) -- cycle;
    \draw[->] (-.5,0) -- (2.5,0);
    \draw (1,-1) node {Case c};
    \draw[->] (0,-1) -- (0,1);
    \draw (.5,-.1) -- (.5,.1);
    \draw (0,1) node[left]{$x$};
    \draw (2.5,0) node[right]{$t$};
    \draw (.5,0) node[below]{$t_n$};
    \draw (2,-.1) -- (2,.1);
    \draw (2,0) node[above]{$t_{n+1}$};
  \end{tikzpicture}
  \]
  \caption{Three cases for discrete transitions. Exact trajectories are depicted
    in dark gray, the over-approximated flow pipes in light gray.}
  \label{fig:transitions-possibilities}
\end{figure}
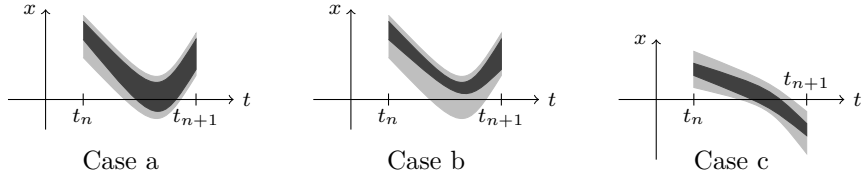

The function $\gsolveZC{}$ is described in
Algorithm~\ref{alg:guaranteed-zc-algorithm} and runs as
follows. First, if for all edges $e\in l_n^\bullet$, $G(e)(\hat
x_n^h)=\bot$, then no transition was activated between $t_n$ and
$t_{n+1}$, and we do nothing (lines 2--4). Otherwise, if there is
$e\in l_n^\bullet$ that may have been activated, then we make sure
that we have $G(e)(\hat x_{n+})=\top$, \ie that the event really
occurred between $\hat t_n$ and $\hat t_{n+1}$ (this is the case (c)
in Figure~\ref{fig:transitions-possibilities}, other cases are handled
as ``special cases'' below), which is achieved by continuing the
guaranteed integration of $F(l_n)$ until we have $G(e)(\hat
x_{n+1})$. This is the role of the while loop (lines 6--10), in which
we also compute the hull of all Picard over-approximations computed
during this process. Then, we are sure that $e$ occurred between~$\hat
x_n$ and~$\hat x_{n+1}$. We then reduce the time interval $[\hat
t_n,\hat t_{n+1}]$ in order to precisely enclose the time $\hat
t_{zc}$ at which the condition $G(e)$ became true (line~11). To do so,
we use the guaranteed polynomial extrapolation~$p$ of
Section~\ref{sec:guaranteed-polynomial-interpolation} to approximate
the value of $x$ between $\hat t_n$ and $\hat t_{n+1}$ without having
to call $\gsolveODE{}$, and use a bisection algorithm to find the
lower and upper limits of $\hat t_{zc}$.

To get the lower limit (the upper limit is obtained similarly), the
bisection algorithm perform as follows. We start with a working list
containing $[\hat t_n,\hat t_{n+1}]$, the convex hull of both time
stamps. Then, we pick the first element $\hat t$ of the working list
and evaluate $p$ on it. If $p(\hat t)=\botop$ and the width of $\hat
t$ is larger than the desired precision, we split $\hat t$ into $\hat
t_1$ and $\hat t_2$ and add them to the working list. If the width
$\hat t$ is smaller than the precision, we return $\hat t$. If $p(\hat
t)=\bot$, we discard $\hat t$ and continue with the rest of the
working list. Note that we cannot have $p(\hat t)=\top$.  The method
to find the upper limit is the same, except that we discard $\hat t$
if $p(\hat t)=\top$.

Finally, once we have $\hat t_{zc}$, we use the guaranteed polynomial
again to compute the zero-crossing state $\hat x_{zc}=p(\hat t_{zc})$
and set $\hat x_{n+1}=R(e)(\hat x_{zc})$, \ie{} we apply the reset
map.

Notice that our algorithm needs to maintain the invariant $G(e)(\hat
x_n)=\bot$ for all $e\in l_n^\bullet$. This imposes that we sometimes
have a particular formulation for zero-crossing conditions. For
instance, the guard and reset functions of the windy ball of
example~\ref{ex:windy-ball} should be reformulated as $G(e)=x<0$ and
$R(e)(x,y,v)=(x,0,-0.8v)$. Under this new formulation, just after the
zero-crossing action has been performed, we have $x=0$ and therefore
the zero-crossing condition $x<0$ is not true. Otherwise, with the
first formulation, the simulation will fail at first
zero-crossing. The transformation is performed automatically for usual
conditions in HySon. Note also that it may be the case that there
exist $e'\in l_{n+1}^\bullet$ such that $G(e')(\hat x_{n+1})\neq
\bot$, \ie{} a transition starting from $l_{n+1}$ may be activated by
$\hat x_{n+1}$. In this case, we execute the transition immediately
after $e$, and continue until we arrive in a location $l$ such that no
transition starting from $l$ is activated. We assume that such $l$
exists, which is true if the \HA $\mc{H}$ does not have Zeno behavior.

\begin{algorithm}[t]
  \caption{Guaranteed Zero-crossing algorithm}
    \label{alg:guaranteed-zc-algorithm}
  \begin{algorithmic}[1]
    \Require{$\mc{H}=(L,F,E,G,R)$, a hybrid automaton}
    \Function{$\gsolveZC{}$}{$\hat x_n,\hat x_{n+1},\hat x^h_n,t_n,h_n,l_n$}
    \If{$\forall e\in l_n^\bullet,\ G(e)(\hat x_n^h)=\bot$}
\State    \Return{$\hat x_{n+1},l_n,t_n+h_n$} \Comment{No jumps}
    \EndIf
    \State Let $e=(l_n,l_{n+1})\in l_n^\bullet$ be such that $G(e)(\hat x_n^h)=\botop$
    \While{$G(e)(\hat x_{n+1})\neq \top$}
    \Let{$(\hat x_{n+1},\hat x^h)$}{$\gsolveODE{}(F(l_n),\hat
        x_{n+1},h_n)$}
      \Let{$x_n^h$}{$x_n^h\cup x^h$}
      \Let{$h_n$}{$h_n+h_n$}
      \EndWhile\Comment{Now $G(e)(\hat x_n)=\bot$ and $G(e)(\hat
        x_{n+1})=\top$}
      \Let{$t_{zc}$}{$\mathrm{tightInterval}(\hat x_n,\hat
        x_{n+1},t_n,t_{n+1})$}
      \Let{$x_{zc}$}{$\mathrm{GPolyODE}(\hat x_n,\hat x_{n+1},\hat
        x_n^h,t_{zc})$}
      \State      \Return ($R(e)(\hat x_{zc}),l_{n+1},t_{zc}$)
\EndFunction
  \end{algorithmic}
\end{algorithm}

\bigskip\noindent\textbf{Special Cases.}
If there is more than one transition activated during the step from
$t_k$ to $t_{k+1}$, we first reject the step and continue with a
reduced step-size. This way, we shall eventually reach a step-size
where only one condition is activated and not the other. If we cannot
separate both transitions before reaching a minimal step-size, we use
our previous algorithm on both transitions separately, apply both
reset maps and then we follow both possible trajectories, \ie we have
a disjunctive analysis when we are not sure of the location.

Finally, we shall discuss the case when the state at times $t_k$ and
$t_{k+1}$ do not verify the guard of a transition $e$ but the hull
computed by Picard iteration does (see
Figure~\ref{fig:transitions-possibilities}, cases a and b). Then,
either the trajectories between $t_k$ and $t_{k+1}$ cross twice the
guard boundary and we missed a zero-crossing, (case a) or it is the
over-approximation due to Picard iteration which makes the guard
validated (case b). We use again our bisection algorithm to
distinguish between these two cases and perform a disjunctive analysis
if we cannot differentiate between them.

%% file: experimentation.tex
We implemented our method in a tool named HySon. It is written in OCaml and
takes as input a representation of a hybrid system either using a set of
equations similar to the ones defined in~\cite{BC12} or a \simulink{} model (for
now without stateflow support). We first present the output of HySon on some
continuous or hybrid systems, and then we compare the performances of HySon with
other tools.

\subsection{Continuous Systems}

\begin{figure}[tbp]
  \centering
  \[
  \includegraphics{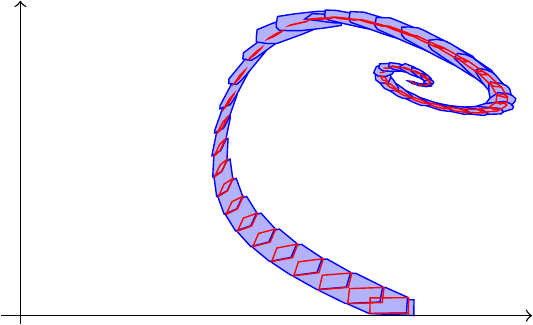}
  \qquad\qquad
  \includegraphics{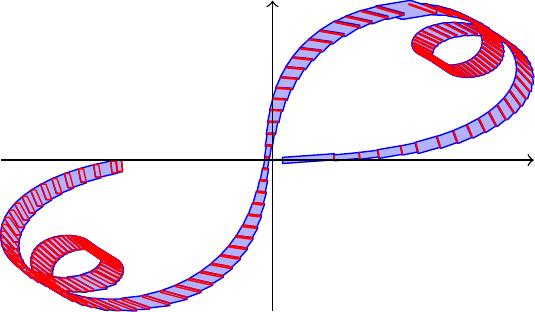}
  \]
  \caption{Over-approximation of the trajectories of the Brusselator
    (left) and Car (right) systems. The blue sets are the
    over-approximations for all $t$ (given by Picard iteration) and the
    red sets are the tight enclosures at the discretization time stamps.}
  \label{fig:tikz-examples-continuous}
\end{figure}

\noindent\textbf{Brusselator.}
We consider the following system, also used in~\cite{CAS12}:
\[
\dot{x}=1 + x^2y - 2.5x
\qquad\qquad
\dot{y}=1.5x - x^2y
\qquad\qquad
x(0)\in[0.9,1]
\qquad\qquad
y(0)\in[0,0.1]
\]
HySon computes the flowpipe up to $t=15$ in $14.3$s, see
Figure~\ref{fig:tikz-examples-continuous}, left.

\bigskip

\noindent {\bf Car.}
We consider the initial value problem given by:
\[
\begin{aligned}
  \dot{x} &= v\cos(0.2t)\cos(\theta)\qquad\qquad
  \dot{y} = v\cos(0.2t)\sin(\theta)\\
  \dot{\theta} &= v\sin(0.2t) / 5
\end{aligned}
\qquad\qquad
\begin{aligned}
  &x(0)=0\qquad y(0)=0\\&\theta(0)=[0,0.1]
\end{aligned}
\]
HySon computes the flowpipe up to $t=30$ in $55.9$s, see
Figure~\ref{fig:tikz-examples-continuous}, right.

\subsection{Hybrid Systems}

We now present two hybrid systems: a ball bouncing on a sinusoidal floor
and a non-linear system with a polynomial jump condition.

\bigskip
\noindent{\bf Ball bouncing on a sinusoidal floor.}
A ball is falling on a sinusoidal floor, and we consider a dynamics with
non-linear wind friction for the ball. The dynamics of the system is
given by
\[
\dot v_x=0
\qquad\qquad
\dot x=v_x
\qquad\qquad
\dot v_y=-g+k v_y^2
\qquad\qquad
\dot y=v_y
\]
starting from the initial conditions $x(0)=1.6$, $v_x(0)=0$, $y(0)=5$
and $v_y(0)=-5$. The bouncing of the ball is given by the transition:
{\small\[
  \pa{
    \begin{aligned}
      v_x&=e(v_d-v_x)\\
      v_y&=e(v_d\cos(x)-v_y)\\
      y&=\sin(x)\\
    \end{aligned}
  }
  \qquad
  \text{when $y<\sin(x)$}
  \]}%
with $v_d=(v_x+v_y\cos(x))/(1+\cos(x)^2)$, where $g=9.8$, $k=0.3$
and $e=0.8$. Note that the exact dynamics of this system is almost
chaotic. HySon is able to compute flow-pipe for this system, as shown on
the following figure.
\[
\begin{tikzpicture}[axis/.style={->,color=black},point/.style={thick},yscale=.6,xscale=1.8,color=red]
  \draw[color=white,pattern=north west lines,pattern color=darkgray] plot[domain=0:6] (\x,{sin(\x r)}) -- (6,-1) -- (0,-1) -- cycle;
  \draw[color=darkgray] plot[domain=0:6] (\x,{sin(\x r)});
  \draw[axis] (0,0) -- (6,0) node[right] {$x$};
  \draw[axis] (0,-0.992226) -- (0,4.949886) node[right] {$y$};
  \input{tikz/ball-sine.tex}
  \input{tikz/ball-sine-matlab.tex}
\end{tikzpicture}
\]

\noindent\textbf{Wolfgram.}
We study the following system, with $a=2$:
\[
\dot{x}(t)=
\begin{cases}
  t^2 + 2x&\text{if }(x + 3/20)^2 + (t + 1/20)^2<1\\
  2t^2 + 3x^2 - a&\text{otherwise}
\end{cases}
\qquad\qquad
x(0)\in[0.3,0.31]
\]
The dynamics of the system is relatively simple, however the jump
condition is a polynomial and is thus not well suited for classical
intersection techniques as in~\cite{CAS12,spaceex}. Our bisection
algorithm for computing the zero-crossing time encloses precisely the
jumping time. To precisely enclose the value of $x$, we insert a reset
in the discrete transition and set $x=\sqrt{1-(t+1/20)^2}-3/20$. This
transformation allows us to obtain a tight enclosure of $x$ as
well. Note however that we performed this transformation manually for
now except for polynomial guard, our future work will include the
automatization of this task for more expressions.

\subsection{Comparison with other Tools}
We now compare the performance of HySon with other tools for reachability
analysis of non-linear hybrid systems: Flow$^*$ as in~\cite{CAS12} and
HydLogic~\cite{hydlogic}. We downloaded both tools from the web and run them on
various examples included in the Flow$^*$ distribution (we could not compile
HydLogic). We run HySon on the same examples and present the execution time for
both in Table~\ref{tab:comparison-hyson-flow}. We see that HySon outperforms
Flow$^*$ on all these examples, whether they are purely continuous systems
(VanDerPol, Brusselator or Lorenz) or hybrid systems (Watertank). Note that for
the Lorenz system, we set a fixed step-size of $0.02$ to achieve a good
precision, which explains the large computation time. For all other examples, we
used a variable step-size and an order $3$ for the Taylor models used in
Flow$^*$.  Let us remark however that some examples work well on Flow$^*$ but
not in HySon, especially the examples with many transitions that may happen
simultaneously. We also want to point out that our tool performs well on linear
examples. We compared it with SpaceEx~\cite{spaceex} on simple examples where
HySon and SpaceEx produced very similar results in terms of precision and
computation time (Appendix~\ref{sec:comparison-spaceex}).

\begin{table}[tbp]
  \centering
  \caption{Experimental results. LOC is the number
    of locations, VAR the number of variables and T the final time of
    simulation. TT is the computation time, in seconds.}
  \label{tab:comparison-hyson-flow}
  \begin{tabular}{|c||c|c|c|c|c|}
    \hline
    Benchmark    &LOC&VAR&T&TT (HySon)&TT (Flow$^*$)\\\hline\hline
    Brusselator&1&2&15&14.3&49.97\\\hline
    Van-der-Pol&1&2&6&16.2&49.17\\\hline
    Lorenz&1&3&1&13.32&119.94\\\hline
    WaterTank&2&5&30&4.35&316.72\\\hline
    Hybrid3D&2&3&2.0&26.65&237.4\\\hline
    Pendulum&1&2&3.8&26.75&N/A\\\hline
    Diode oscillator~\cite{Frehse05}&3&2&20&29.56&42.65\\\hline
  \end{tabular}
\end{table}

%% file: tikz/ball-sine.tex
\draw[point] (1.600000,4.949886) rectangle (1.600000,4.949886);
\draw[point] (1.600000,4.899549) rectangle (1.600000,4.899549);
\draw[point] (1.600000,4.848995) rectangle (1.600000,4.848995);
\draw[point] (1.600000,4.798232) rectangle (1.600000,4.798232);
\draw[point] (1.600000,4.747264) rectangle (1.600000,4.747264);
\draw[point] (1.600000,4.696099) rectangle (1.600000,4.696099);
\draw[point] (1.600000,4.644742) rectangle (1.600000,4.644742);
\draw[point] (1.600000,4.593200) rectangle (1.600000,4.593200);
\draw[point] (1.600000,4.541477) rectangle (1.600000,4.541477);
\draw[point] (1.600000,4.489579) rectangle (1.600000,4.489579);
\draw[point] (1.600000,4.437512) rectangle (1.600000,4.437512);
\draw[point] (1.600000,4.385281) rectangle (1.600000,4.385281);
\draw[point] (1.600000,4.332891) rectangle (1.600000,4.332891);
\draw[point] (1.600000,4.280347) rectangle (1.600000,4.280347);
\draw[point] (1.600000,4.227653) rectangle (1.600000,4.227653);
\draw[point] (1.600000,4.174815) rectangle (1.600000,4.174815);
\draw[point] (1.600000,4.121836) rectangle (1.600000,4.121836);
\draw[point] (1.600000,4.068722) rectangle (1.600000,4.068722);
\draw[point] (1.600000,4.015476) rectangle (1.600000,4.015476);
\draw[point] (1.600000,3.962103) rectangle (1.600000,3.962103);
\draw[point] (1.600000,3.908606) rectangle (1.600000,3.908606);
\draw[point] (1.600000,3.854990) rectangle (1.600000,3.854990);
\draw[point] (1.600000,3.801258) rectangle (1.600000,3.801258);
\draw[point] (1.600000,3.747414) rectangle (1.600000,3.747414);
\draw[point] (1.600000,3.693461) rectangle (1.600000,3.693461);
\draw[point] (1.600000,3.639404) rectangle (1.600000,3.639404);
\draw[point] (1.600000,3.585245) rectangle (1.600000,3.585245);
\draw[point] (1.600000,3.530987) rectangle (1.600000,3.530987);
\draw[point] (1.600000,3.476634) rectangle (1.600000,3.476634);
\draw[point] (1.600000,3.422189) rectangle (1.600000,3.422189);
\draw[point] (1.600000,3.367655) rectangle (1.600000,3.367655);
\draw[point] (1.600000,3.313034) rectangle (1.600000,3.313034);
\draw[point] (1.600000,3.258330) rectangle (1.600000,3.258330);
\draw[point] (1.600000,3.203545) rectangle (1.600000,3.203545);
\draw[point] (1.600000,3.148681) rectangle (1.600000,3.148681);
\draw[point] (1.600000,3.093742) rectangle (1.600000,3.093742);
\draw[point] (1.600000,3.038729) rectangle (1.600000,3.038729);
\draw[point] (1.600000,2.983646) rectangle (1.600000,2.983646);
\draw[point] (1.600000,2.928494) rectangle (1.600000,2.928494);
\draw[point] (1.600000,2.873276) rectangle (1.600000,2.873276);
\draw[point] (1.600000,2.817993) rectangle (1.600000,2.817993);
\draw[point] (1.600000,2.762648) rectangle (1.600000,2.762648);
\draw[point] (1.600000,2.707243) rectangle (1.600000,2.707243);
\draw[point] (1.600000,2.651781) rectangle (1.600000,2.651781);
\draw[point] (1.600000,2.596261) rectangle (1.600000,2.596261);
\draw[point] (1.600000,2.540688) rectangle (1.600000,2.540688);
\draw[point] (1.600000,2.485062) rectangle (1.600000,2.485062);
\draw[point] (1.600000,2.429385) rectangle (1.600000,2.429385);
\draw[point] (1.600000,2.373659) rectangle (1.600000,2.373659);
\draw[point] (1.600000,2.317885) rectangle (1.600000,2.317885);
\draw[point] (1.600000,2.262065) rectangle (1.600000,2.262065);
\draw[point] (1.600000,2.206201) rectangle (1.600000,2.206201);
\draw[point] (1.600000,2.150293) rectangle (1.600000,2.150293);
\draw[point] (1.600000,2.094345) rectangle (1.600000,2.094345);
\draw[point] (1.600000,2.038355) rectangle (1.600000,2.038355);
\draw[point] (1.600000,1.982327) rectangle (1.600000,1.982327);
\draw[point] (1.600000,1.926262) rectangle (1.600000,1.926262);
\draw[point] (1.600000,1.870160) rectangle (1.600000,1.870160);
\draw[point] (1.600000,1.814022) rectangle (1.600000,1.814022);
\draw[point] (1.600000,1.757851) rectangle (1.600000,1.757851);
\draw[point] (1.600000,1.701647) rectangle (1.600000,1.701647);
\draw[point] (1.600000,1.645411) rectangle (1.600000,1.645411);
\draw[point] (1.600000,1.589145) rectangle (1.600000,1.589145);
\draw[point] (1.600000,1.532848) rectangle (1.600000,1.532848);
\draw[point] (1.600000,1.476523) rectangle (1.600000,1.476523);
\draw[point] (1.600000,1.420170) rectangle (1.600000,1.420170);
\draw[point] (1.600000,1.363791) rectangle (1.600000,1.363791);
\draw[point] (1.600000,1.307385) rectangle (1.600000,1.307385);
\draw[point] (1.600000,1.250954) rectangle (1.600000,1.250954);
\draw[point] (1.600000,1.194499) rectangle (1.600000,1.194499);
\draw[point] (1.600000,1.138021) rectangle (1.600000,1.138021);
\draw[point] (1.600000,1.081520) rectangle (1.600000,1.081520);
\draw[point] (1.600000,1.024997) rectangle (1.600000,1.024997);
\draw[point] (1.600000,0.999574) rectangle (1.600000,0.999574);
\draw[point] (1.602639,1.044545) rectangle (1.602639,1.044545);
\draw[point] (1.605279,1.089139) rectangle (1.605279,1.089139);
\draw[point] (1.607918,1.133343) rectangle (1.607918,1.133343);
\draw[point] (1.610557,1.177149) rectangle (1.610557,1.177149);
\draw[point] (1.613197,1.220545) rectangle (1.613197,1.220545);
\draw[point] (1.615836,1.263521) rectangle (1.615836,1.263521);
\draw[point] (1.618476,1.306065) rectangle (1.618476,1.306065);
\draw[point] (1.621115,1.348166) rectangle (1.621115,1.348166);
\draw[point] (1.623754,1.389814) rectangle (1.623754,1.389814);
\draw[point] (1.626394,1.430996) rectangle (1.626394,1.430996);
\draw[point] (1.629033,1.471701) rectangle (1.629033,1.471701);
\draw[point] (1.631672,1.511918) rectangle (1.631672,1.511918);
\draw[point] (1.634312,1.551633) rectangle (1.634312,1.551633);
\draw[point] (1.636951,1.590836) rectangle (1.636951,1.590836);
\draw[point] (1.639590,1.629514) rectangle (1.639591,1.629514);
\draw[point] (1.642230,1.667654) rectangle (1.642230,1.667654);
\draw[point] (1.644869,1.705244) rectangle (1.644869,1.705244);
\draw[point] (1.647509,1.742272) rectangle (1.647509,1.742272);
\draw[point] (1.650148,1.778725) rectangle (1.650148,1.778725);
\draw[point] (1.652787,1.814591) rectangle (1.652787,1.814591);
\draw[point] (1.655427,1.849856) rectangle (1.655427,1.849856);
\draw[point] (1.658066,1.884507) rectangle (1.658066,1.884507);
\draw[point] (1.660705,1.918532) rectangle (1.660706,1.918532);
\draw[point] (1.663345,1.951918) rectangle (1.663345,1.951919);
\draw[point] (1.665984,1.984653) rectangle (1.665984,1.984653);
\draw[point] (1.668623,2.016722) rectangle (1.668624,2.016722);
\draw[point] (1.671263,2.048113) rectangle (1.671263,2.048113);
\draw[point] (1.673902,2.078813) rectangle (1.673902,2.078813);
\draw[point] (1.676542,2.108810) rectangle (1.676542,2.108810);
\draw[point] (1.679181,2.138090) rectangle (1.679181,2.138090);
\draw[point] (1.681820,2.166641) rectangle (1.681821,2.166642);
\draw[point] (1.684460,2.194451) rectangle (1.684460,2.194451);
\draw[point] (1.687099,2.221506) rectangle (1.687099,2.221506);
\draw[point] (1.689738,2.247795) rectangle (1.689739,2.247795);
\draw[point] (1.692378,2.273305) rectangle (1.692378,2.273305);
\draw[point] (1.695017,2.298024) rectangle (1.695017,2.298024);
\draw[point] (1.697657,2.321941) rectangle (1.697657,2.321941);
\draw[point] (1.700296,2.345043) rectangle (1.700296,2.345044);
\draw[point] (1.702935,2.367321) rectangle (1.702936,2.367321);
\draw[point] (1.705575,2.388761) rectangle (1.705575,2.388761);
\draw[point] (1.708214,2.409354) rectangle (1.708214,2.409354);
\draw[point] (1.710853,2.429089) rectangle (1.710854,2.429089);
\draw[point] (1.713493,2.447956) rectangle (1.713493,2.447956);
\draw[point] (1.716132,2.465945) rectangle (1.716132,2.465945);
\draw[point] (1.718771,2.483046) rectangle (1.718772,2.483046);
\draw[point] (1.721411,2.499250) rectangle (1.721411,2.499251);
\draw[point] (1.724050,2.514549) rectangle (1.724051,2.514550);
\draw[point] (1.726690,2.528934) rectangle (1.726690,2.528935);
\draw[point] (1.729329,2.542398) rectangle (1.729329,2.542398);
\draw[point] (1.731968,2.554932) rectangle (1.731969,2.554932);
\draw[point] (1.734608,2.566529) rectangle (1.734608,2.566530);
\draw[point] (1.737247,2.577184) rectangle (1.737247,2.577184);
\draw[point] (1.739886,2.586890) rectangle (1.739887,2.586890);
\draw[point] (1.742526,2.595642) rectangle (1.742526,2.595642);
\draw[point] (1.745165,2.603434) rectangle (1.745166,2.603434);
\draw[point] (1.747804,2.610263) rectangle (1.747805,2.610263);
\draw[point] (1.750444,2.616123) rectangle (1.750444,2.616123);
\draw[point] (1.753083,2.621012) rectangle (1.753084,2.621012);
\draw[point] (1.755208,2.624241) rectangle (1.755209,2.624241);
\draw[point] (1.756713,2.626144) rectangle (1.756713,2.626145);
\draw[point] (1.758033,2.627552) rectangle (1.758033,2.627553);
\draw[point] (1.759740,2.629011) rectangle (1.759741,2.629012);
\draw[point] (1.762004,2.630314) rectangle (1.762004,2.630315);
\draw[point] (1.764643,2.630922) rectangle (1.764644,2.630924);
\draw[point] (1.766389,2.630786) rectangle (1.766390,2.630789);
\draw[point] (1.767780,2.630371) rectangle (1.767781,2.630374);
\draw[point] (1.770420,2.628835) rectangle (1.770420,2.628840);
\draw[point] (1.773059,2.626320) rectangle (1.773060,2.626326);
\draw[point] (1.775218,2.623536) rectangle (1.775218,2.623544);
\draw[point] (1.777857,2.619246) rectangle (1.777858,2.619255);
\draw[point] (1.780497,2.613983) rectangle (1.780497,2.613994);
\draw[point] (1.782696,2.608857) rectangle (1.782696,2.608869);
\draw[point] (1.785335,2.601818) rectangle (1.785336,2.601832);
\draw[point] (1.787974,2.593816) rectangle (1.787975,2.593832);
\draw[point] (1.790614,2.584856) rectangle (1.790614,2.584874);
\draw[point] (1.793253,2.574942) rectangle (1.793254,2.574962);
\draw[point] (1.795892,2.564081) rectangle (1.795893,2.564103);
\draw[point] (1.798532,2.552279) rectangle (1.798532,2.552302);
\draw[point] (1.801171,2.539542) rectangle (1.801172,2.539567);
\draw[point] (1.803811,2.525877) rectangle (1.803811,2.525904);
\draw[point] (1.806450,2.511292) rectangle (1.806451,2.511320);
\draw[point] (1.809089,2.495795) rectangle (1.809090,2.495825);
\draw[point] (1.811729,2.479394) rectangle (1.811729,2.479426);
\draw[point] (1.814368,2.462098) rectangle (1.814369,2.462132);
\draw[point] (1.817007,2.443917) rectangle (1.817008,2.443953);
\draw[point] (1.819647,2.424860) rectangle (1.819647,2.424898);
\draw[point] (1.822286,2.404937) rectangle (1.822287,2.404976);
\draw[point] (1.824925,2.384159) rectangle (1.824926,2.384199);
\draw[point] (1.827565,2.362535) rectangle (1.827566,2.362577);
\draw[point] (1.830204,2.340077) rectangle (1.830205,2.340120);
\draw[point] (1.832844,2.316796) rectangle (1.832844,2.316841);
\draw[point] (1.835483,2.292703) rectangle (1.835484,2.292750);
\draw[point] (1.838122,2.267810) rectangle (1.838123,2.267858);
\draw[point] (1.840762,2.242129) rectangle (1.840762,2.242179);
\draw[point] (1.843401,2.215672) rectangle (1.843402,2.215724);
\draw[point] (1.846040,2.188452) rectangle (1.846041,2.188504);
\draw[point] (1.848680,2.160479) rectangle (1.848681,2.160533);
\draw[point] (1.851319,2.131768) rectangle (1.851320,2.131824);
\draw[point] (1.853959,2.102331) rectangle (1.853959,2.102387);
\draw[point] (1.856598,2.072179) rectangle (1.856599,2.072237);
\draw[point] (1.859237,2.041327) rectangle (1.859238,2.041387);
\draw[point] (1.861877,2.009787) rectangle (1.861877,2.009848);
\draw[point] (1.864516,1.977572) rectangle (1.864517,1.977634);
\draw[point] (1.867155,1.944695) rectangle (1.867156,1.944758);
\draw[point] (1.869795,1.911168) rectangle (1.869796,1.911233);
\draw[point] (1.872434,1.877005) rectangle (1.872435,1.877071);
\draw[point] (1.875073,1.842219) rectangle (1.875074,1.842286);
\draw[point] (1.877713,1.806822) rectangle (1.877714,1.806890);
\draw[point] (1.880352,1.770827) rectangle (1.880353,1.770896);
\draw[point] (1.882992,1.734248) rectangle (1.882992,1.734318);
\draw[point] (1.885631,1.697096) rectangle (1.885632,1.697167);
\draw[point] (1.888270,1.659385) rectangle (1.888271,1.659457);
\draw[point] (1.890910,1.621126) rectangle (1.890911,1.621200);
\draw[point] (1.893549,1.582333) rectangle (1.893550,1.582407);
\draw[point] (1.896188,1.543017) rectangle (1.896189,1.543093);
\draw[point] (1.898828,1.503192) rectangle (1.898829,1.503268);
\draw[point] (1.901467,1.462868) rectangle (1.901468,1.462945);
\draw[point] (1.904106,1.422057) rectangle (1.904107,1.422136);
\draw[point] (1.906746,1.380773) rectangle (1.906747,1.380852);
\draw[point] (1.909385,1.339025) rectangle (1.909386,1.339105);
\draw[point] (1.912025,1.296826) rectangle (1.912025,1.296907);
\draw[point] (1.914664,1.254187) rectangle (1.914665,1.254269);
\draw[point] (1.917303,1.211119) rectangle (1.917304,1.211201);
\draw[point] (1.919943,1.167632) rectangle (1.919944,1.167716);
\draw[point] (1.922582,1.123739) rectangle (1.922583,1.123823);
\draw[point] (1.925221,1.079449) rectangle (1.925222,1.079533);
\draw[point] (1.927861,1.034772) rectangle (1.927862,1.034857);
\draw[point] (1.930500,0.989719) rectangle (1.930501,0.989805);
\draw[point] (1.933139,0.944300) rectangle (1.933140,0.944387);
\draw[point] (1.933684,0.934873) rectangle (1.933690,0.934876);
\draw[point] (1.958360,0.961519) rectangle (1.958367,0.961522);
\draw[point] (1.983035,0.987391) rectangle (1.983043,0.987396);
\draw[point] (2.007711,1.012478) rectangle (2.007720,1.012484);
\draw[point] (2.032386,1.036768) rectangle (2.032396,1.036775);
\draw[point] (2.057062,1.060250) rectangle (2.057073,1.060258);
\draw[point] (2.081737,1.082911) rectangle (2.081749,1.082920);
\draw[point] (2.106413,1.104740) rectangle (2.106426,1.104750);
\draw[point] (2.131088,1.125727) rectangle (2.131102,1.125739);
\draw[point] (2.155764,1.145861) rectangle (2.155779,1.145874);
\draw[point] (2.180439,1.165131) rectangle (2.180455,1.165145);
\draw[point] (2.205115,1.183528) rectangle (2.205132,1.183543);
\draw[point] (2.229790,1.201042) rectangle (2.229808,1.201058);
\draw[point] (2.254466,1.217663) rectangle (2.254485,1.217681);
\draw[point] (2.279141,1.233382) rectangle (2.279161,1.233402);
\draw[point] (2.303817,1.248192) rectangle (2.303838,1.248212);
\draw[point] (2.328492,1.262083) rectangle (2.328514,1.262105);
\draw[point] (2.353168,1.275048) rectangle (2.353191,1.275072);
\draw[point] (2.377843,1.287081) rectangle (2.377867,1.287106);
\draw[point] (2.402519,1.298173) rectangle (2.402544,1.298199);
\draw[point] (2.427194,1.308320) rectangle (2.427220,1.308347);
\draw[point] (2.451870,1.317514) rectangle (2.451897,1.317543);
\draw[point] (2.476545,1.325752) rectangle (2.476573,1.325782);
\draw[point] (2.501221,1.333027) rectangle (2.501250,1.333059);
\draw[point] (2.525896,1.339336) rectangle (2.525926,1.339369);
\draw[point] (2.550572,1.344676) rectangle (2.550602,1.344710);
\draw[point] (2.575248,1.349043) rectangle (2.575279,1.349078);
\draw[point] (2.599923,1.352434) rectangle (2.599955,1.352471);
\draw[point] (2.624599,1.354847) rectangle (2.624632,1.354887);
\draw[point] (2.642602,1.355991) rectangle (2.642636,1.356031);
\draw[point] (2.666222,1.356699) rectangle (2.666257,1.356742);
\draw[point] (2.674919,1.356734) rectangle (2.674955,1.356778);
\draw[point] (2.686276,1.356596) rectangle (2.686312,1.356642);
\draw[point] (2.699566,1.356171) rectangle (2.699603,1.356218);
\draw[point] (2.724242,1.354629) rectangle (2.724279,1.354678);
\draw[point] (2.748917,1.352107) rectangle (2.748956,1.352160);
\draw[point] (2.768928,1.349345) rectangle (2.768967,1.349400);
\draw[point] (2.793604,1.345055) rectangle (2.793644,1.345113);
\draw[point] (2.818279,1.339792) rectangle (2.818320,1.339853);
\draw[point] (2.839102,1.334596) rectangle (2.839144,1.334660);
\draw[point] (2.863778,1.327547) rectangle (2.863820,1.327614);
\draw[point] (2.888453,1.319535) rectangle (2.888497,1.319606);
\draw[point] (2.913129,1.310565) rectangle (2.913173,1.310639);
\draw[point] (2.937804,1.300641) rectangle (2.937850,1.300719);
\draw[point] (2.962480,1.289770) rectangle (2.962526,1.289851);
\draw[point] (2.987155,1.277958) rectangle (2.987203,1.278042);
\draw[point] (3.011831,1.265211) rectangle (3.011879,1.265298);
\draw[point] (3.036506,1.251537) rectangle (3.036556,1.251627);
\draw[point] (3.061182,1.236942) rectangle (3.061232,1.237036);
\draw[point] (3.085857,1.221436) rectangle (3.085909,1.221532);
\draw[point] (3.110533,1.205026) rectangle (3.110585,1.205125);
\draw[point] (3.135208,1.187721) rectangle (3.135262,1.187823);
\draw[point] (3.159884,1.169530) rectangle (3.159938,1.169636);
\draw[point] (3.184559,1.150464) rectangle (3.184615,1.150573);
\draw[point] (3.209235,1.130532) rectangle (3.209291,1.130644);
\draw[point] (3.233910,1.109745) rectangle (3.233968,1.109859);
\draw[point] (3.258586,1.088112) rectangle (3.258644,1.088230);
\draw[point] (3.283261,1.065645) rectangle (3.283321,1.065766);
\draw[point] (3.307937,1.042356) rectangle (3.307997,1.042479);
\draw[point] (3.332612,1.018255) rectangle (3.332674,1.018381);
\draw[point] (3.357288,0.993354) rectangle (3.357350,0.993482);
\draw[point] (3.381963,0.967665) rectangle (3.382027,0.967796);
\draw[point] (3.406639,0.941200) rectangle (3.406703,0.941334);
\draw[point] (3.431314,0.913971) rectangle (3.431380,0.914108);
\draw[point] (3.455990,0.885991) rectangle (3.456056,0.886130);
\draw[point] (3.480665,0.857272) rectangle (3.480733,0.857414);
\draw[point] (3.505341,0.827827) rectangle (3.505409,0.827971);
\draw[point] (3.530016,0.797668) rectangle (3.530086,0.797815);
\draw[point] (3.554692,0.766809) rectangle (3.554762,0.766958);
\draw[point] (3.579367,0.735262) rectangle (3.579438,0.735413);
\draw[point] (3.604043,0.703040) rectangle (3.604115,0.703193);
\draw[point] (3.628718,0.670155) rectangle (3.628791,0.670311);
\draw[point] (3.653394,0.636622) rectangle (3.653468,0.636780);
\draw[point] (3.678070,0.602453) rectangle (3.678144,0.602613);
\draw[point] (3.702745,0.567660) rectangle (3.702821,0.567822);
\draw[point] (3.727421,0.532257) rectangle (3.727497,0.532421);
\draw[point] (3.752096,0.496256) rectangle (3.752174,0.496422);
\draw[point] (3.776772,0.459670) rectangle (3.776850,0.459839);
\draw[point] (3.801447,0.422513) rectangle (3.801527,0.422683);
\draw[point] (3.826123,0.384796) rectangle (3.826203,0.384968);
\draw[point] (3.850798,0.346531) rectangle (3.850880,0.346706);
\draw[point] (3.875474,0.307733) rectangle (3.875556,0.307909);
\draw[point] (3.900149,0.268412) rectangle (3.900233,0.268590);
\draw[point] (3.924825,0.228581) rectangle (3.924909,0.228760);
\draw[point] (3.949500,0.188252) rectangle (3.949586,0.188433);
\draw[point] (3.974176,0.147437) rectangle (3.974262,0.147619);
\draw[point] (3.998851,0.106147) rectangle (3.998939,0.106331);
\draw[point] (4.023527,0.064395) rectangle (4.023615,0.064581);
\draw[point] (4.048202,0.022191) rectangle (4.048292,0.022378);
\draw[point] (4.072878,-0.020453) rectangle (4.072968,-0.020264);
\draw[point] (4.097553,-0.063525) rectangle (4.097645,-0.063335);
\draw[point] (4.122229,-0.107016) rectangle (4.122321,-0.106824);
\draw[point] (4.146904,-0.150914) rectangle (4.146998,-0.150721);
\draw[point] (4.171580,-0.195208) rectangle (4.171674,-0.195014);
\draw[point] (4.196255,-0.239889) rectangle (4.196351,-0.239693);
\draw[point] (4.220931,-0.284946) rectangle (4.221027,-0.284748);
\draw[point] (4.245606,-0.330368) rectangle (4.245704,-0.330170);
\draw[point] (4.270282,-0.376147) rectangle (4.270380,-0.375947);
\draw[point] (4.294957,-0.422272) rectangle (4.295057,-0.422072);
\draw[point] (4.319633,-0.468735) rectangle (4.319733,-0.468533);
\draw[point] (4.344308,-0.515525) rectangle (4.344410,-0.515322);
\draw[point] (4.368984,-0.562634) rectangle (4.369086,-0.562430);
\draw[point] (4.393659,-0.610052) rectangle (4.393763,-0.609847);
\draw[point] (4.418335,-0.657772) rectangle (4.418439,-0.657566);
\draw[point] (4.443010,-0.705784) rectangle (4.443115,-0.705577);
\draw[point] (4.467686,-0.754081) rectangle (4.467792,-0.753873);
\draw[point] (4.492361,-0.802654) rectangle (4.492468,-0.802445);
\draw[point] (4.517037,-0.851495) rectangle (4.517145,-0.851285);
\draw[point] (4.541712,-0.900597) rectangle (4.541821,-0.900386);
\draw[point] (4.566388,-0.949951) rectangle (4.566498,-0.949740);
\draw[point] (4.587388,-0.992226) rectangle (4.587619,-0.992198);
\draw[point] (4.616283,-0.958853) rectangle (4.616500,-0.958835);
\draw[point] (4.645178,-0.926132) rectangle (4.645381,-0.926126);
\draw[point] (4.674074,-0.894087) rectangle (4.674262,-0.894070);
\draw[point] (4.702969,-0.862723) rectangle (4.703143,-0.862691);
\draw[point] (4.731864,-0.832051) rectangle (4.732024,-0.832002);
\draw[point] (4.760759,-0.802083) rectangle (4.760906,-0.802017);
\draw[point] (4.789654,-0.772831) rectangle (4.789787,-0.772749);
\draw[point] (4.818549,-0.744309) rectangle (4.818668,-0.744211);
\draw[point] (4.847444,-0.716530) rectangle (4.847549,-0.716414);
\draw[point] (4.876339,-0.689505) rectangle (4.876430,-0.689371);
\draw[point] (4.905234,-0.663247) rectangle (4.905312,-0.663096);
\draw[point] (4.934130,-0.637769) rectangle (4.934193,-0.637599);
\draw[point] (4.963025,-0.613081) rectangle (4.963074,-0.612893);
\draw[point] (4.991920,-0.589197) rectangle (4.991955,-0.588990);
\draw[point] (5.020815,-0.566127) rectangle (5.020836,-0.565901);
\draw[point] (5.049698,-0.543883) rectangle (5.049730,-0.543638);
\draw[point] (5.078578,-0.522476) rectangle (5.078626,-0.522212);
\draw[point] (5.107458,-0.501917) rectangle (5.107522,-0.501633);
\draw[point] (5.136338,-0.482217) rectangle (5.136419,-0.481913);
\draw[point] (5.165217,-0.463385) rectangle (5.165315,-0.463061);
\draw[point] (5.194097,-0.445431) rectangle (5.194212,-0.445087);
\draw[point] (5.222977,-0.428366) rectangle (5.223108,-0.428001);
\draw[point] (5.251857,-0.412197) rectangle (5.252004,-0.411812);
\draw[point] (5.280737,-0.396935) rectangle (5.280901,-0.396528);
\draw[point] (5.309617,-0.382587) rectangle (5.309797,-0.382159);
\draw[point] (5.338497,-0.369160) rectangle (5.338693,-0.368711);
\draw[point] (5.367377,-0.356664) rectangle (5.367590,-0.356193);
\draw[point] (5.396257,-0.345103) rectangle (5.396486,-0.344611);
\draw[point] (5.425137,-0.334486) rectangle (5.425382,-0.333972);
\draw[point] (5.454017,-0.324818) rectangle (5.454279,-0.324282);
\draw[point] (5.482897,-0.316105) rectangle (5.483175,-0.315547);
\draw[point] (5.511777,-0.308351) rectangle (5.512071,-0.307771);
\draw[point] (5.540657,-0.301562) rectangle (5.540968,-0.300959);
\draw[point] (5.569537,-0.295740) rectangle (5.569864,-0.295115);
\draw[point] (5.598417,-0.290890) rectangle (5.598760,-0.290242);
\draw[point] (5.623954,-0.287412) rectangle (5.624312,-0.286745);
\draw[point] (5.652668,-0.284414) rectangle (5.653043,-0.283724);
\draw[point] (5.679021,-0.282514) rectangle (5.679410,-0.281803);
\draw[point] (5.705920,-0.281415) rectangle (5.706324,-0.280682);
\draw[point] (5.720360,-0.281175) rectangle (5.720773,-0.280431);
\draw[point] (5.734934,-0.281182) rectangle (5.735355,-0.280426);
\draw[point] (5.745342,-0.281340) rectangle (5.745769,-0.280575);
\draw[point] (5.756383,-0.281646) rectangle (5.756816,-0.280873);
\draw[point] (5.785263,-0.283125) rectangle (5.785712,-0.282328);
\draw[point] (5.814142,-0.285582) rectangle (5.814609,-0.284761);
\draw[point] (5.842431,-0.288937) rectangle (5.842913,-0.288093);
\draw[point] (5.871311,-0.293327) rectangle (5.871809,-0.292459);
\draw[point] (5.900191,-0.298690) rectangle (5.900706,-0.297799);
\draw[point] (5.919881,-0.302903) rectangle (5.920408,-0.301995);
\draw[point] (5.948761,-0.309895) rectangle (5.949304,-0.308963);
\draw[point] (5.962878,-0.314858) rectangle (5.965808,-0.312076);
\draw[point] (5.957692,-0.291868) rectangle (5.960651,-0.289068);
\draw[point] (5.952506,-0.269705) rectangle (5.955494,-0.266886);
\draw[point] (5.947320,-0.248380) rectangle (5.950337,-0.245542);
\draw[point] (5.942134,-0.227904) rectangle (5.945180,-0.225047);
\draw[point] (5.936948,-0.208287) rectangle (5.940023,-0.205412);
\draw[point] (5.931762,-0.189540) rectangle (5.934866,-0.186645);
\draw[point] (5.926576,-0.171672) rectangle (5.929709,-0.168758);
\draw[point] (5.921390,-0.154693) rectangle (5.924552,-0.151759);
\draw[point] (5.916204,-0.138612) rectangle (5.919395,-0.135658);
\draw[point] (5.911017,-0.123437) rectangle (5.914238,-0.120464);
\draw[point] (5.905831,-0.109178) rectangle (5.909081,-0.106185);
\draw[point] (5.900645,-0.095842) rectangle (5.903924,-0.092828);
\draw[point] (5.895459,-0.083435) rectangle (5.898767,-0.080401);
\draw[point] (5.890273,-0.071966) rectangle (5.893610,-0.068911);
\draw[point] (5.885087,-0.061441) rectangle (5.888453,-0.058365);
\draw[point] (5.879901,-0.051865) rectangle (5.883296,-0.048768);
\draw[point] (5.874715,-0.043245) rectangle (5.878139,-0.040127);
\draw[point] (5.869529,-0.035584) rectangle (5.872982,-0.032445);
\draw[point] (5.864343,-0.028888) rectangle (5.867825,-0.025728);
\draw[point] (5.859156,-0.023160) rectangle (5.862668,-0.019979);
\draw[point] (5.853970,-0.018405) rectangle (5.857511,-0.015201);
\draw[point] (5.848948,-0.014728) rectangle (5.852517,-0.011504);
\draw[point] (5.843762,-0.011892) rectangle (5.847360,-0.008646);
\draw[point] (5.839471,-0.010285) rectangle (5.843093,-0.007022);
\draw[point] (5.835040,-0.009330) rectangle (5.838687,-0.006047);
\draw[point] (5.832447,-0.009103) rectangle (5.836109,-0.005809);
\draw[point] (5.830015,-0.009112) rectangle (5.833691,-0.005808);
\draw[point] (5.828169,-0.009263) rectangle (5.831854,-0.005951);
\draw[point] (5.825889,-0.009621) rectangle (5.829588,-0.006299);
\draw[point] (5.820703,-0.011141) rectangle (5.824431,-0.007796);
\draw[point] (5.815517,-0.013639) rectangle (5.819274,-0.010272);
\draw[point] (5.810619,-0.016896) rectangle (5.814404,-0.013508);
\draw[point] (5.805433,-0.021293) rectangle (5.809247,-0.017882);
\draw[point] (5.800247,-0.026662) rectangle (5.804090,-0.023229);
\draw[point] (5.796265,-0.031444) rectangle (5.800130,-0.027993);
\draw[point] (5.791079,-0.038525) rectangle (5.794973,-0.035052);
\draw[point] (5.785893,-0.046569) rectangle (5.789816,-0.043073);
\draw[point] (5.780706,-0.055572) rectangle (5.784659,-0.052053);
\draw[point] (5.775520,-0.065527) rectangle (5.779502,-0.061986);
\draw[point] (5.770334,-0.076430) rectangle (5.774345,-0.072867);
\draw[point] (5.765148,-0.088273) rectangle (5.769188,-0.084689);
\draw[point] (5.759962,-0.101052) rectangle (5.764031,-0.097445);
\draw[point] (5.754776,-0.114757) rectangle (5.758874,-0.111129);
\draw[point] (5.749590,-0.129383) rectangle (5.753717,-0.125733);
\draw[point] (5.744404,-0.144920) rectangle (5.748560,-0.141249);
\draw[point] (5.739218,-0.161360) rectangle (5.743403,-0.157668);
\draw[point] (5.734032,-0.178695) rectangle (5.738246,-0.174982);
\draw[point] (5.728845,-0.196915) rectangle (5.733089,-0.193182);
\draw[point] (5.723659,-0.216010) rectangle (5.727932,-0.212257);
\draw[point] (5.718473,-0.235971) rectangle (5.722775,-0.232198);
\draw[point] (5.713287,-0.256788) rectangle (5.717618,-0.252994);
\draw[point] (5.708101,-0.278448) rectangle (5.712461,-0.274635);
\draw[point] (5.702915,-0.300943) rectangle (5.707304,-0.297110);
\draw[point] (5.697729,-0.324260) rectangle (5.702147,-0.320408);
\draw[point] (5.692543,-0.348388) rectangle (5.696990,-0.344518);
\draw[point] (5.687357,-0.373316) rectangle (5.691833,-0.369427);
\draw[point] (5.682171,-0.399031) rectangle (5.686676,-0.395124);
\draw[point] (5.676984,-0.425522) rectangle (5.681519,-0.421597);
\draw[point] (5.671798,-0.452777) rectangle (5.676362,-0.448833);
\draw[point] (5.666612,-0.480782) rectangle (5.671205,-0.476821);
\draw[point] (5.661426,-0.509525) rectangle (5.666048,-0.505547);
\draw[point] (5.656240,-0.538995) rectangle (5.660891,-0.535000);
\draw[point] (5.651054,-0.569177) rectangle (5.655734,-0.565165);
\draw[point] (5.645345,-0.595471) rectangle (5.653061,-0.589245);
\draw[point] (5.620072,-0.594692) rectangle (5.627970,-0.588505);
\draw[point] (5.601584,-0.594742) rectangle (5.609616,-0.588585);
\draw[point] (5.589032,-0.595076) rectangle (5.597154,-0.588937);
\draw[point] (5.572927,-0.595857) rectangle (5.581165,-0.589744);
\draw[point] (5.547654,-0.597886) rectangle (5.556074,-0.591811);
\draw[point] (5.522381,-0.600892) rectangle (5.530983,-0.594856);
\draw[point] (5.501730,-0.604074) rectangle (5.510481,-0.598069);
\draw[point] (5.476457,-0.608853) rectangle (5.485390,-0.602886);
\draw[point] (5.451183,-0.614603) rectangle (5.460299,-0.608675);
\draw[point] (5.425910,-0.621322) rectangle (5.435208,-0.615431);
\draw[point] (5.400637,-0.629005) rectangle (5.410117,-0.623152);
\draw[point] (5.375364,-0.637649) rectangle (5.385026,-0.631833);
\draw[point] (5.350091,-0.647247) rectangle (5.359935,-0.641468);
\draw[point] (5.324818,-0.657795) rectangle (5.334844,-0.652053);
\draw[point] (5.299544,-0.669286) rectangle (5.309753,-0.663581);
\draw[point] (5.274271,-0.681714) rectangle (5.284662,-0.676046);
\draw[point] (5.248998,-0.695073) rectangle (5.259571,-0.689441);
\draw[point] (5.223725,-0.709354) rectangle (5.234480,-0.703758);
\draw[point] (5.198452,-0.724550) rectangle (5.209389,-0.718989);
\draw[point] (5.173179,-0.740652) rectangle (5.184298,-0.735127);
\draw[point] (5.147905,-0.757652) rectangle (5.159207,-0.752162);
\draw[point] (5.122632,-0.775541) rectangle (5.134116,-0.770086);
\draw[point] (5.097359,-0.794309) rectangle (5.109025,-0.788888);
\draw[point] (5.072086,-0.813946) rectangle (5.083934,-0.808559);
\draw[point] (5.046813,-0.834443) rectangle (5.058843,-0.829089);
\draw[point] (5.021540,-0.855788) rectangle (5.033752,-0.850467);
\draw[point] (4.996267,-0.877971) rectangle (5.008661,-0.872682);
\draw[point] (4.970993,-0.900981) rectangle (4.983570,-0.895724);
\draw[point] (4.945720,-0.924806) rectangle (4.958479,-0.919581);
\draw[point] (4.920447,-0.949435) rectangle (4.933388,-0.944241);
\draw[point] (4.895174,-0.974856) rectangle (4.908297,-0.969693);
\draw[point] (4.878812,-0.986326) rectangle (4.901475,-0.982161);
\draw[point] (4.852827,-0.973666) rectangle (4.875246,-0.970222);
\draw[point] (4.826842,-0.961942) rectangle (4.849017,-0.959223);
\draw[point] (4.800857,-0.951160) rectangle (4.822788,-0.949171);
\draw[point] (4.774872,-0.941902) rectangle (4.796559,-0.939496);
\draw[point] (4.748887,-0.934185) rectangle (4.770330,-0.930191);
\draw[point] (4.722901,-0.927433) rectangle (4.744101,-0.921843);
\draw[point] (4.696916,-0.921649) rectangle (4.717872,-0.914456);
\draw[point] (4.673297,-0.917235) rectangle (4.694031,-0.908579);
\draw[point] (4.647312,-0.913309) rectangle (4.667802,-0.903039);
\draw[point] (4.621326,-0.910359) rectangle (4.641573,-0.898471);
\draw[point] (4.595341,-0.908389) rectangle (4.615344,-0.894878);
\draw[point] (4.569740,-0.907405) rectangle (4.589503,-0.892294);
\draw[point] (4.543755,-0.907380) rectangle (4.563274,-0.890642);
\draw[point] (4.523892,-0.908022) rectangle (4.543224,-0.890040);
\draw[point] (4.506419,-0.909060) rectangle (4.525587,-0.889984);
\draw[point] (4.484074,-0.911032) rectangle (4.503033,-0.890557);
\draw[point] (4.458273,-0.914210) rectangle (4.476990,-0.892120);
\draw[point] (4.432288,-0.918382) rectangle (4.450761,-0.894670);
\draw[point] (4.406303,-0.923529) rectangle (4.424532,-0.898197);
\draw[point] (4.380317,-0.929646) rectangle (4.398302,-0.902699);
\draw[point] (4.300439,-0.944163) rectangle (4.372681,-0.916245);
\draw[point] (4.263752,-0.945890) rectangle (4.372927,-0.900834);

%% file: tikz/ball-sine-matlab.tex
\draw [color=blue,opacity=0.5] (1.60000000000000,5) -- (1.60000000000000,4.74483021964011) -- (1.60000000000000,4.23132062889773) -- (1.60000000000000,3.69974745186396) -- (1.60000000000000,3.15507173565258) -- (1.60000000000000,2.60272710292816) -- (1.60000000000000,2.04487506797174) -- (1.60000000000000,1.48308131357397) -- (1.60000000000000,0.999573603041505) -- (1.60000000000000,0.999573603041465) -- (1.60000000000000,0.999573603041506) -- (1.62639404849163,1.43099668993056) -- (1.65278809698326,1.81459171553910) -- (1.67918214547488,2.13809186647493) -- (1.70557619396651,2.38876345408018) -- (1.73197024245814,2.55493487431747) -- (1.75836429094977,2.62786972535061) -- (1.78475833944139,2.60345762947193) -- (1.81115238793302,2.48309388558329) -- (1.83754643642465,2.27337486830975) -- (1.86394048491628,1.98474160661266) -- (1.89033453340790,1.62961829083328) -- (1.91672858189953,1.22066218490624) -- (1.93369291881194,0.934872504173727) -- (1.93369291881194,0.934872504173662) -- (1.93369291881198,0.934872504173700) -- (2.17728860512552,1.16270972744430) -- (2.42088429143906,1.30581460662881) -- (2.66764855202833,1.35673049004309) -- (2.91441281261759,1.31012423774465) -- (3.16117707320685,1.16864725297522) -- (3.40794133379612,0.939920276228265) -- (3.65470559438538,0.635005031112780) -- (3.90146985497464,0.266518765614042) -- (4.14823411556391,-0.153024794754118) -- (4.38971335250711,-0.602146145464815) -- (4.58758903016250,-0.992222588527778) -- (4.58758903016264,-0.992222588528061) -- (4.58758903016272,-0.992222588527965) -- (4.86550584517097,-0.699511343693567) -- (5.14342266017923,-0.477254844237186) -- (5.43225455493416,-0.331562488636948) -- (5.72108644968910,-0.280490369511529) -- (5.96458233743720,-0.313240143478108) -- (5.96458233743728,-0.313240143478131) -- (5.96458233743726,-0.313240143478032) -- (5.91278025740047,-0.121766574045138) -- (5.86097817736368,-0.0214321894783957) -- (5.80917609732689,-0.0178108597244591) -- (5.75737401729010,-0.111109918602793) -- (5.70557193725331,-0.296124796896382) -- (5.65376985721652,-0.563336781660113) -- (5.64883465574180,-0.592654617186267) -- (5.64883465574177,-0.592654617186444) -- (5.64883465574151,-0.592654617186431) -- (5.39688746315303,-0.628995909261697) -- (5.14494027056455,-0.760741952423621) -- (4.89299307797607,-0.980736713350721) -- (4.88950261180050,-0.984356339148088) -- (4.88950261180036,-0.984356339148237) -- (4.88950261180012,-0.984356339148120) -- (4.62849990743441,-0.903938969994598) -- (4.36749720306870,-0.920741623123006) -- (4.33503417281517,-0.929642541499016) -- (4.33503417281501,-0.929642541499062) -- (4.33503417281494,-0.929642541498951) -- (4.33503417281487,-0.929642541498847) -- (4.33503417281480,-0.929642541498735) -- (4.33503417281473,-0.929642541498623) -- (4.33503417281458,-0.929642541498399) -- (4.33503417281429,-0.929642541497951) -- (4.33503417281371,-0.929642541497055) -- (4.33503417281255,-0.929642541495264) -- (4.33503417281022,-0.929642541491681) -- (4.33503417280557,-0.929642541484516) -- (4.33503417279626,-0.929642541470185) -- (4.33503417277765,-0.929642541441523) -- (4.33503417274043,-0.929642541384199) -- (4.33503417266599,-0.929642541269552) -- (4.33503417251710,-0.929642541040257) -- (4.33503417221932,-0.929642540581667) -- (4.33503417162377,-0.929642539664486) -- (4.33503417043266,-0.929642537830126) -- (4.33503416805045,-0.929642534161405) -- (4.33503416328603,-0.929642526823963) -- (4.33503415375718,-0.929642512149081) -- (4.33503413469949,-0.929642482799317) -- (4.33503409658411,-0.929642424099797) -- (4.33503402035336,-0.929642306700783) -- (4.33503386789184,-0.929642071902864) -- (4.33503356296880,-0.929641602307463) -- (4.33503295312273,-0.929640663118404) -- (4.33503173343058,-0.929638784747263) -- (4.33502929404628,-0.929635028032886) -- (4.33502441527769,-0.929627514715751) -- (4.33501465774050,-0.929612488527964) -- (4.33499514266613,-0.929582437938337) -- (4.33495611251738,-0.929522343903023) -- (4.33487805221988,-0.929402184409356) -- (4.33472193162490,-0.929161979739481) -- (4.33440969043492,-0.928682027746442) -- (4.33378520805496,-0.927723953761656) -- (4.33253624329504,-0.925815130703727) -- (4.33003831377521,-0.922026823340335) -- (4.32504245473555,-0.914567874115739);

%% file: conclusion.tex
We presented a new approach to compute the flowpipes of nonlinear hybrid
systems using guaranteed version of numerical methods. Our method is
based on guaranteed explicit Runge-Kutta integration methods and on a
new guaranteed polynomial interpolation based on the well-known
Hermite-Birkoff method. This interpolation is cheap and precise to
over-approximate continuous state values. Using both methods, we can
precisely compute flowpipes of nonlinear hybrid systems, with a few
number of restrictions on the nature of flows and jumps. Remark that
with guaranteed polynomial interpolation, we can accurately and soundly
handle nonlinear jumps in hybrid systems without using an intersection
operator which is usually costly to define. Note also that we can handle
in the same manner invariants in hybrid automaton using our algorithm
for zero-crossing events. More precisely, we would add a new step in the
simulation loop to check that the invariant is fulfilled at each
integration step. Finally, the experiments showed that our approach is
efficient and precise on a set of representative case studies: we showed
that our approach outperforms existing techniques on the flowpipe
computation of nonlinear systems.

As future work, we plan to handle multiple zero-crossing events
involving trajectories associated to different system behaviors. As a
result, to keep the flowpipe computation sharp we must handle
disjunctive futures efficiently. We also want to extend our parser of
Simulink models, presented in~\cite{BCM12}, to handle Stateflow and thus
apply our tool on more realistic examples.